\documentclass[preprint,12pt]{elsarticle} 
\usepackage[english]{babel}

\usepackage{enumerate}           
\usepackage{amsmath}
\usepackage{amssymb}    
\usepackage{graphicx}      
\usepackage{epsfig}       
\usepackage{algorithm} 
\usepackage{algpseudocode}
\usepackage{booktabs}   
\usepackage{subcaption}

\usepackage[most]{tcolorbox}

\newtcolorbox{highlightblock}{          
  enhanced,
  colback=yellow!25,
  colframe=yellow!60!black,
  boxrule=0pt,
  sharp corners,
  breakable,       
  left=2pt, right=2pt, top=2pt, bottom=2pt
}

\usepackage{color}      
\usepackage{caption} 
\usepackage{float}
\usepackage{tabulary} 
\usepackage{csquotes}

\usepackage{enumerate}   

\usepackage[a4paper,top=2.0cm,bottom=2.0cm,left=1.0cm,right=1.0cm]{geometry}

\journal{Int. J. Numer. Meth. Eng.}
\begin{document}

\title{Piecewise Linear Approximation and PID Control Optimization for Nonlinear Systems}
\author{Robert~Vrabel}
\ead{robert.vrabel@stuba.sk}
\address{Slovak University of Technology in Bratislava, Institute of Applied Informatics, Automation and Mechatronics,  Bottova 25,  917 24 Trnava, Slovakia\\ (ORCID: 0000-0002-2640-595X)}

\begin{abstract}
	This paper investigates the control of nonlinear systems using a piecewise linear approximation framework. The proposed approach combines a PID controller with locally linearized models obtained by partitioning the nonlinear function into subregions over a compact domain. This approximation yields an analytically tractable representation of the system dynamics, enabling the application of transfer-based and frequency-domain analysis tools that are not directly applicable to nonlinear systems. 
	
	As the number of linear segments increases, the approximated system progressively approaches the behavior of the original nonlinear system, allowing for a meaningful frequency-domain interpretation of the dynamics. The PID controller parameters are optimized using the Particle Swarm Optimization method with performance criteria based on ITAE (Integral of Time-weighted Absolute Error) and ISO (Integral of Squared Overshoot). 
	
	Numerical simulations confirm the effectiveness of the proposed method, demonstrating that controller parameters obtained from the piecewise linear model ensure stable and accurate control when applied to the original nonlinear system, while maintaining a balance between computational effort and approximation accuracy.
\end{abstract}

              
\newtheorem{thm}{Theorem} 
\newtheorem{lem}[thm]{Lemma}
\newtheorem{defn}[thm]{Definition}
\newdefinition{rem}{Remark}
\newdefinition{ex}{Example}
\newproof{pf}{Proof}
\newproof{pot1}{Proof of Theorem \ref{thm1}} 
\newproof{pot2}{Proof of Theorem \ref{thm2}}

\maketitle

\section{Introduction}
Transfer functions are a fundamental tool in classical control theory and linear system analysis, providing a frequency-domain representation that greatly simplifies the study of system behavior. In linear systems, they allow for systematic design and tuning of controllers, such as PID controllers, to achieve desired transient and steady-state performance. However, for nonlinear systems, the concept of a transfer function cannot be directly applied due to the system's inherent nonlinearities. Despite this limitation, several methodologies have been developed to extend the transfer function concept to nonlinear systems \cite{schetzen1980volterra}, enabling approximate linear analysis and controller design within a nonlinear framework.

The describing function method approximates certain nonlinear systems by linearizing them around a specific operating point, particularly for systems with sinusoidal inputs. This approach provides an approximate frequency-domain analysis but is limited to specific types of nonlinearities and input signals. 

The Volterra series generalizes the concept of transfer functions to nonlinear systems by representing the system's output as a series of multidimensional convolutions of the input. Each term in the series corresponds to a higher-order frequency response function, capturing the system's nonlinear behavior. This method is particularly useful for weakly nonlinear systems. 

Relatively recent research has introduced the concept of quasi-linear transfer functions, which aim to characterize the output frequency behavior of nonlinear systems. This method provides a new perspective on analyzing nonlinear systems in the frequency domain. 

Another approach involves extending system-level synthesis methods to polynomial dynamical systems, yielding finite impulse response, time-invariant, closed-loop transfer functions with guaranteed disturbance cancellation. This generalizes feedback linearization to enable partial feedback linearization over a finite-time horizon. 

These resources offer a comprehensive overview of the methodologies developed to analyze nonlinear systems using concepts analogous to transfer functions. 

The proportional-integral-derivative (PID) controller is a cornerstone in feedback control systems due to its simplicity, effectiveness, and versatility. It is a general-purpose algorithm used to regulate processes in a wide variety of applications, from industrial automation to robotics, avionics, and beyond. The controller works by continuously adjusting the control input to minimize the error between a desired setpoint and the actual process output. The optimization of PID control parameters for nonlinear systems is more challenging than for linear systems because nonlinearities can cause unpredictable behavior, such as oscillations (limit cycles, chaos-induced oscillations, self-sustained oscillations), instability, or slow convergence. Nonlinear systems often exhibit behavior such as saturation, dead-zone, hysteresis, and time-varying dynamics that complicate the control process. Traditional PID controllers are designed for linear systems, and their performance in nonlinear systems may degrade due to the inability to model nonlinearities directly. Thus, optimal PID tuning in nonlinear systems is crucial for avoiding issues like overshoot, oscillations, or slow response. To achieve optimal PID performance in nonlinear systems, several methods and strategies are employed.

In some cases, the nonlinear system can be linearized around an operating point. This simplified linear model can then be used to apply classical PID tuning methods. However, this method only works for relatively small deviations from the operating point and does not always capture the full nonlinear behavior (if the equilibrium point is not hyperbolic, for example).

Adaptive PID controllers adjust the PID parameters in real-time based on the changing dynamics of the nonlinear system. This approach is useful when the nonlinear system's behavior changes with time or operating conditions. Adaptive methods can be based on model reference adaptive control (MRAC), where the controller parameters are modified to match a desired reference model in real-time to account for changing system dynamics or operating conditions. The paper \cite{Zuo2017adaptivePID} proposes an adaptive piecewise linear switch controller that combines PID and MRAC to enhance the response time and tracking performance of hydraulic actuators, which are inherently nonlinear systems.

Gain scheduling involves adjusting the PID parameters based on the operating conditions of the nonlinear system. This approach requires identifying different operating regimes of the system, each of which may need different PID gains. For example, the PID gains can be varied depending on system speed, temperature, or load.

Optimization algorithms such as Genetic Algorithms, Particle Swarm Optimization (PSO), and Simulated Annealing can be used to find the optimal PID parameters for nonlinear systems by minimizing a cost function, typically involving criteria like settling time, overshoot, steady-state error and control effort.

In some cases, using a pure PID controller may not be effective for nonlinear systems. In such cases, nonlinear control methods such as fuzzy logic controllers, sliding mode controllers, or backstepping control are often considered. However, PID controllers can still be used in combination with these methods to improve performance. By incorporating fuzzy logic into the PID structure, the controller can handle nonlinearities more effectively. The fuzzy logic component allows the controller to adjust its behavior based on the input-output relationship, compensating for nonlinearities.

Neural networks can also be used to tune the PID parameters or modify the control strategy. This technique can handle complex nonlinearities by learning from data.
 
Book \cite{astrom2006advanced} discusses PID control design and optimization for various systems, including nonlinear systems, with practical insights into tuning and implementation.  A widely-used textbook \cite{nise2011control} offers a solid foundation for both linear and nonlinear control systems, including techniques for PID control of nonlinear systems. Book \cite{ogata2010modern} provides a detailed treatment of control system analysis and design, including PID control and nonlinear system behavior. 

The Particle Swarm Optimization (PSO) algorithm as initially developed by James Kennedy and Russell Eberhart \cite{eberhart1995particle,kennedy1997discrete,shi1998modified}, inspired by the collective movement of bird flocks and fish schools. The idea was to simulate the social behavior of these animals to solve optimization problems. The initial paper \cite{eberhart1995particle} presented the basic PSO model, which was intended to simulate the social behavior of organisms in nature and solve nonlinear optimization problems. PSO's history is characterized by its rapid evolution and widespread adoption in diverse fields. Initially designed for simple optimization problems, it has since become a robust tool for solving complex, multi-dimensional, and dynamic optimization challenges. Researchers continue to explore new variants and applications, ensuring PSO's relevance in modern computational intelligence. Over time, various PSO extensions have been developed to address multi-objective optimization problems \cite{coello2002mopso,deb2013evolutionary}. These multi-objective PSO methods have since been widely adopted in real-world engineering and scientific optimization tasks. Several studies started combining PSO with other optimization algorithms (such as genetic algorithms or simulated annealing) to create hybrid models that could perform better on certain problems. PSO was adapted for constrained optimization problems, where researchers developed specific mechanisms to handle constraints in the search space \cite{poli2007particle}. PSO has gained widespread use in machine learning for tasks such as feature selection, training deep learning networks, and hyperparameter tuning. PSO has been applied in various industries such as robotics, control systems, image processing, and power systems optimization. In response to large-scale optimization problems, parallel and distributed PSO algorithms have been developed to take advantage of modern computing architectures \cite{zhang2018comprehensive}. The basic PSO algorithm has been modified in several ways to limit the particle movement and avoid divergence, to control the influence of previous velocity on the next velocity, improving convergence and exploration/exploitation balance and adjusting these parameters to balance local search (individual experience) and global search (group experience). The study \cite{Charkoutsis2023} introduces a novel nonlinear PID (NLPID) controller, optimized using PSO, tailored for first-order plus time delay systems. The proposed controller employs time-varying gains to enhance both set-point tracking and disturbance rejection capabilities. Simulation results indicate that the PSO-tuned NLPID controller outperforms traditional PID controllers, offering faster response times, minimized overshoot, and improved robustness against parametric uncertainties. 

\begin{figure}
\begin{center}
\includegraphics[width=\textwidth]{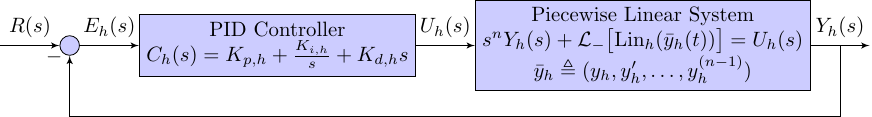}      
\caption{Feedback control system with piecewise linear approximation of the nonlinear system.}  
\label{fig1_control}                                  
\end{center}                                 
\end{figure}

Unlike the aforementioned articles, this paper develops a method applicable to systems of arbitrary order. Conceptually, the theory is mathematically supported; however, in implementation, technical challenges may arise, such as the algorithm for partitioning a compact set in \( \mathbb{R}^n \) into simplices for higher dimensions \( n \). The objective of this paper is to develop a methodology for the analysis and synthesis of control circuits by approximating a nonlinear system in a feedback configuration, as illustrated in the accompanying schematic (Figure~\ref{fig1_control}), using transfer functions derived via piecewise linear, more precisely affine approximation. We demonstrate that, under specific assumptions, the solutions of the sequence of linearized systems are well-defined over the entire real-time axis and uniformly converge to the solution of the original nonlinear model. Similarly, the corresponding transfer functions exhibit analogous convergence behavior, and the limit function can be legitimately interpreted as the transfer function of the nonlinear system, the formula (\ref{transfer:nonlinear}).

\begin{rem}
Although a classical PID controller can, in principle, be directly implemented in a feedback loop with a nonlinear plant, such an approach lacks a rigorous analytical foundation for determining controller parameters or guaranteeing closed-loop stability. For nonlinear systems, the superposition principle and frequency-domain concepts used in conventional control theory no longer apply, which prevents the direct derivation of transfer functions and stability margins. Consequently, the tuning of controller parameters in the nonlinear case typically relies on empirical or heuristic optimization methods, which may be computationally demanding and provide limited theoretical insight.

To address this issue, the nonlinear system is approximated by a locally linearized (piecewise linear) model, which enables the use of \emph{transfer algebra} and frequency-domain synthesis tools that are not directly applicable to nonlinear dynamics. This approximation provides an analytical framework for controller design, facilitates the derivation of closed-form stability conditions, and significantly reduces computational complexity. 
Furthermore, for the piecewise linearized control system, it becomes possible to compute the \emph{initial jump} and \emph{initial slope} at \(t = 0^+\) exactly, using the procedure derived in Remark~\ref{initialJumpandSlope} and Appendix~\ref{app:IJandIS}. These quantities characterize the instantaneous and short-term transient response of the closed-loop system and are of particular importance in practical controller implementation, as they allow accurate prediction of the initial system behavior immediately after actuation, ensuring smooth and stable control performance. 
\end{rem}

In the subsequent stage of the study, we address the design of PID controller parameters (\( K_p, K_i, K_d \)) through a cost function minimization approach utilizing PSO algorithm, implemented via the {\tt pyswarm} module which provides a flexible and efficient way to use PSO in Python. The efficacy of the proposed methodology is validated by comparing its performance against traditional controller parameter tuning methods, such as closed-loop tuning, in the context of a linear time-invariant system. Furthermore, in Subsection~\ref{example_nonlinear}, we illustrate the application of the proposed approach to a nonlinear system, showing that the PID parameters optimized for its piecewise linear approximation yield a comparable response when applied to the original nonlinear model.

First, let us introduce two important technical results that will be used in the proofs of existence and uniqueness of the solutions, as well as their uniform convergence.
\section[]{Banach Contraction Principle and Gronwall Inequality}
Both the Banach Contraction Principle \cite{Banach1922} and Gronwall Inequality \cite{Gronwall1919} are fundamental tools in analysis, particularly in the study of differential equations, stability analysis, and control theory. They are often applied to prove existence and uniqueness of solutions, estimate bounds on solutions, and establish stability results.
\begin{thm}[Banach Contraction Principle] 
Let \((X, d)\) be a complete metric space, and let \(T: X \mapsto X\) be a contraction mapping, meaning there exists a constant 
\(c \in [0, 1)\) such that
\[
d(T(x), T(y)) \leq c \cdot d(x, y), \quad \text{for all} \quad x, y \in X.
\]
Then, there exists a unique point \(x^* \in X\) such that \(T(x^*) = x^*\).
\end{thm}
The Gronwall inequality is a fundamental result used to derive bounds on solutions of differential inequalities. It provides a systematic way to estimate the growth of a solution to a differential equation or an inequality, especially when the solution is constrained by an upper bound. There are several forms of Gronwall inequality, but one of the most common forms relevant to differential equations is the differential inequality form. A detailed proof of the version of Gronwall inequality used in this paper can be found in Appendix~\ref{app:Gronwall}.
\begin{thm}[Gronwall Inequality]\label{Gronwall}
Suppose \(w\) is differentiable on an interval \([a,b]\) and satisfies
\[
\frac{dw(t)}{dt} \le f(t) w(t) + g(t), \qquad t\in[a,b],
\]
where \(f\) and \(g\) are continuous on \([a,b]\). Then for every \(t\in[a,b]\)
\[
w(t) \le w(a)\exp\!\Big(\int_a^t f(s)\,ds\Big)
+ \int_a^t g(s)\exp\!\Big(\int_s^t f(\tau)\,d\tau\Big)\,ds .
\]
\end{thm}
\section[]{Linear interpolation of a multivariate function over an $n$-simplex}
To establish the uniform convergence of piecewise-linear approximations for a \( C^2 \) function \( f : \mathbb{R}^n \to \mathbb{R} \), we show that the sequence 
of piecewise-linear interpolants converges uniformly to \( f \) as the partition of the domain is progressively refined. In general, the following analysis applies to twice continuously differentiable functions defined on \( \mathbb{R}^n \); however, whenever possible, the regularity assumptions on \( f \) are relaxed in specific steps to extend the applicability of the results.

\medskip
{\bf Step 1: Define the Piecewise Linear Approximation}

Let \( D \subset \mathbb{R}^n \) be a compact domain, and assume that \( f \) is differentiable on \( D \). We define a sequence of piecewise linear approximations to \( f \), denoted as \( \text{Lin}_h \), constructed as follows:

\medskip
{\bf{(a)}} Partition the domain: We partition the domain \( D \) into a grid of small regions, typically line segments (in \( \mathbb{R}^1 \)), triangles (in \( \mathbb{R}^2 \)), tetrahedrons (in \( \mathbb{R}^3 \)) or, in general, $n-$simplices (in \( \mathbb{R}^n \)), such that the maximum diameter of the cells in the partition becomes arbitrarily small as \( h \) increases, $\{h \}$ is an infinite and increasing sequence of natural numbers. Denote this partition as \( P_h \), where the size of the cells $C_{h,i}$, $i=1,2,\dots,h$ in the partition decreases as \( h \to \infty \). The $n-$simplex $C_{h,i}$ is defined by its vertices \( \{ {P_1}_{h,i}, {P_2}_{h,i}, \dots, {P_{n+1}}_{h,i} \} \) and the function \( f \) is known at these vertices.

\medskip
{\bf{(b)}} Linear approximation on each region: For each cell \( C_{h,i} \in P_h \), we approximate \( f \) by a linear function \( \text{Lin}_{h,i}(y) \) that takes the values of \( f \) at the vertices of the cell. More precisely, on each cell \( C_{h,i} \), we define \( \text{Lin}_{h,i}(y) \) by interpolating the values of \( f \) at the vertices of $C_{h,i}$.
The goal is to compute the value of the function \( \text{Lin}_{h,i}(y) \) at any point \( y \) within the simplex $C_{h,i}$.

In the case of an \( n \)-simplex, the coordinates \( \lambda_1, \lambda_2, \dots, \lambda_{n+1} \) of a point \( y \) inside the simplex can be expressed as the barycentric coordinates. These coordinates are non-negative and satisfy the condition
\[
\sum_{i=1}^{n+1} \lambda_i = 1
\]
The general form for the linear interpolation of the function \( f \) on $C_{h,i}$ is
\[
\text{Lin}_{h,i}(y) = \sum_{j=1}^{n+1} \lambda_j f({P_j}_{h,i}), \ i=1,2,\dots,h
\]
To compute the barycentric coordinates \( \lambda_1, \lambda_2, \dots, \lambda_{n+1} \) for a point \( y \), we need to calculate the areas (or volumes, depending on the dimension) of the sub-simplices formed by \( y \) and the vertices of the simplex.

For a given point \( y \) inside the \( n \)-simplex \( C_{h,i} \), the barycentric coordinates \( \lambda_j \) are proportional to the volume of the \( n-1 \)-simplex formed by excluding the vertex \( {P_j}_{h,i} \), and this is normalized by the volume of the full \( n \)-simplex.

The formula for the barycentric coordinate \( \lambda_j \) in terms of the volume is
\[
\lambda_j = \frac{\text{Vol}(C_{h,i}^{(j)})}{\text{Vol}(C_{h,i})}
\]
where

- \( \text{Vol}(C_{h,i}) \) is the volume of the \( n \)-simplex formed by the vertices \( {P_1}_{h,i}, {P_2}_{h,i}, \dots, {P_{n+1}}_{h,i} \).

- \( \text{Vol}(C_{h,i}^{(j)}) \) is the volume of the sub-simplex formed by the point \( y \) and all the vertices of the simplex except \( {P_j}_{h,i} \).

Thus,
\begin{equation}\label{Lin_h}
\text{Lin}_{h}(y)\triangleq\{ \text{Lin}_{h,i}(y)\, | \,  i=1,2,\dots,h \}
\end{equation}
represents piecewise linear approximation of the function $f$ on the compact domain $D$.

\medskip
{\bf Step 2: Uniform convergence condition}

To prove that \( \text{Lin}_h \) converges uniformly to \( f \), we must show that
\[
\lim_{h \to \infty} \sup_{{y} \in D} |f({y})-\text{Lin}_{h}(y)| = 0.
\]
In other words, the difference between the piecewise linear approximation and \( f \) must become arbitrarily small uniformly for all \( {y} \in D \) as \( h \to \infty \).

The key to proving uniform convergence is to bound the error between \( \text{Lin}_h \) and \( f \) on each cell of the partition and show that the maximum error across all cells tends to zero as the mesh is refined.

The error estimate for piecewise linear interpolation of a differentiable function \( f(y_1, y_2, \dots, y_n) \) depends on the smoothness of the function and the size of the partition (mesh) used in the interpolation. Specifically, when we approximate \( f \) using piecewise linear functions, the error can be described in terms of the diameter of the regions in the partition and the second derivatives of \( f \) \cite{burden2010numerical} as
\[
| f({y}) - \text{Lin}_{h,i}({y}) | \leq \frac{1}{2} \| \nabla^2 f(\xi_{h,i}) \| (\text{diam}(C_{h,i}))^2,
\]
where \( \text{diam}(C_{h,i}) \) is the diameter (the maximum distance between any two points) of the cell \( C_{h,i} \) and \( \nabla^2 f(\xi_{h,i}) \) is the Hessian of \( f \) evaluated at the point \( \xi_{h,i} \) inside the cell.
\( P_h \) be a partition of \( D \) into \( h \) regions, and \( \text{Lin}_{h}({y}) \) be the corresponding piecewise linear interpolation. The global error for the entire domain can be bounded by considering the worst-case error over all the cells \( C_{h,i} \) in the partition. Since the error on each cell \( C_{h,i} \) depends on \( (\text{diam}(C_{h,i}))^2 \), we can sum the errors across all cells. If we refine the partition (that is, make the cells smaller), the error in each cell decreases. The global error is controlled by the maximum diameter of the cells in the partition, so the total error can be bounded by
\[
\sup_{{y} \in D} | f({y}) - \text{Lin}_h({y}) |  \leq \frac{1}{2} \| \nabla^2 f(\xi) \| \cdot \max_{C_{h,i} \in P_h} \left( \text{diam}(C_{h,i}) \right)^2.
\]
Thus, the error estimate for piecewise linear interpolation of \( f(y_1, y_2, \dots, y_n) \) can be expressed as
\[
| f({y}) - \text{Lin}_h({y}) | = O\left((\text{diam}(C_{h,i}))^2\right),
\]
where \( \text{diam}(C_{h,i}) \) is the diameter of the cells in the partition and the constant hidden in the \( O(\cdot) \) notation depends on the smoothness of \( f \), particularly the bounds on \( \|\nabla^2 f\| \) on \( D \). If $f$ is twice continuously differentiable, that is, \( f \in C^2(D\mapsto\mathbb{R}^n) \), then Hessian matrix is bounded and symmetric on $D$  and \( \|\nabla^2 f\| \) represents the largest absolute eigenvalue of \( \nabla^2 f \) on \( D \).
\begin{defn}
A function \( \tilde f: \mathbb{R}^n \mapsto \mathbb{R}^n \) is globally Lipschitz continuous on \( \mathbb{R}^n \) if there exists a constant 
$\tilde{L}$ such that
\[
| \tilde f(x) - \tilde f(y) | \leq \tilde{L} | x - y |, \quad \forall x, y \in \mathbb{R}^n.
\]
Here, \( | \cdot | \) denotes the Euclidean norm, and the smallest such \(\tilde L \) is called the Lipschitz constant of \( \tilde f(y) \).
\end{defn}
If the function \( \tilde f: \mathbb{R}^n \mapsto \mathbb{R}^n \) is continuously differentiable (that is, \( \tilde f \in C^1(\mathbb{R}^n) \)), then a sufficient condition for \( \tilde f \) to be globally Lipschitz is that there exists a constant \( M > 0 \) such that the norm of the Jacobian matrix of \( \tilde f \) is bounded, that is,
\[
\|\nabla \tilde f(y)\| \leq M, \quad \forall y \in \mathbb{R}^n,
\]
where \( \nabla \tilde f(y) \) is the Jacobian matrix of \( \tilde f \) at \( y \) and \( \| \nabla \tilde f(y) \| \) is the operator norm of the Jacobian which is defined as the largest singular value of the Jacobian matrix. This singular value corresponds to the square root of the largest eigenvalue of the matrix \( \nabla \tilde{f}(y)^\top \nabla \tilde{f}(y) \), where the superscript \( \top \) denotes the transpose of the Jacobian matrix \( \nabla \tilde{f}(y) \). If this condition is satisfied, then \( f \) is globally Lipschitz with Lipschitz constant \( \tilde L = M \) because of the mean value theorem.

Consider nonlinear system
\begin{equation}\label{nonlinear_system_n_order}
y^{(n)}(t) + f(y(t), y'(t), \dots, y^{(n-1)}(t)) = u(t),
\end{equation}
on the time interval \([0, \infty)\) with the initial conditions
\[
y(0) = y'(0) = \dots = y^{(n-1)}(0) = 0,
\]
or reformulated as an equivalent system of first-order ordinary differential equations (ODEs) by introducing variables
\[
y_1 = y, \ y_2 = y', \ \dots, \ y_{n-1} = y^{(n-2)},\ y_n = y^{(n-1)}.
\]
Then, the equation becomes
\begin{equation}\label{nonlinear_system_1_order}
\begin{aligned}
\left.
\begin{array}{ll}
y_1' &= y_2, \\
y_2' &= y_3, \\
&\ \vdots \\
y_{n-1}' &= y_n, \\
y_n' &= -f(y_1, y_2, \dots, y_n) + u(t)
\end{array}
\right\} \triangleq y'=\tilde f(y)+\tilde u(t).
\end{aligned}  
\end{equation}
The initial conditions become
\[
y_1(0) = y_2(0) = \dots = y_n(0) = 0.
\]
\begin{lem}\label{Lipschitz1}
If $f$ is globally Lipschitz with Lipschitz constant \(  L \) , then also $\tilde f$ is globally Lipschitz with Lipschitz constant  $\tilde L = ({n-1+L^2})^{1/2}$.
\end{lem}
\begin{pf}
The detailed proof of this lemma can be found in Appendix~\ref{app:Lipschitz}.
\end{pf}
\section[]{Lipschitz Constant of the Piecewise Linear Appro\-xi\-ma\-tion $\text{Lin}_h(y)$}
To prove that the Lipschitz constant \( L_{\text{pwl},h}\) of the piecewise linear approximation of a function \( f({y}) : \mathbb{R}^n \mapsto \mathbb{R} \) satisfies \( L_{\text{pwl},h} \leq L \), we will proceed as follows.
A function \( f({y}) \) is globally Lipschitz if there exists a constant \( L \geq 0 \) such that
\[
|f({x}) - f({y})| \leq L |x-y|, \quad \forall x,y\in \mathbb{R}^n.
\]
A piecewise linear approximation of \( f({y}) \) is constructed by partitioning the compact domain \( D\subset \mathbb{R}^n \) into a set of non-overlapping cells \( \{C_{h,i}, i=1,2,\dots,h\} \). Within each region \( C_{h,i} \), the function is approximated as a linear function \( \text{Lin}_{h,i}({y}) \), so
\[
f({y}) \approx \text{Lin}_{h,i}({y}), \quad {y} \in C_{h,i}.
\]
For simplicity, assume that the piecewise linear approximation is exact at the vertices of the simplex and \( \text{Lin}_{h,i}({y}) \) is defined as
   \[
   \text{Lin}_{h,i}({y}) = {a}_{h,i} {y} + b_{h,i},
   \]
where \( {a}_{h,i}\in \mathbb{R}^n \) is the gradient (or slope vector) of the linear function in region \( C_{h,i} \).

For future reference, we introduce cumulative notation analogous to (\ref{Lin_h}),  
\[
a_h \triangleq \{a_{h,i} \, | \, i = 1, 2, \dots, h\}
\]  
and  
\[
b_h \triangleq \{b_{h,i} \, | \, i = 1, 2, \dots, h\}.
\]  
This notation will be used to collectively represent the sequences of coefficients \( \{a_{h,i}\} \) and \( \{b_{h,i}\} \) for \( i = 1, 2, \dots, h \).

The Lipschitz constant \( L_{\text{pwl},h} \) of the piecewise linear approximation is defined as
\[
L_{\text{pwl},h} = \sup_{x,y \in D, x \neq y} \frac{|\text{Lin}_{h,i}({x}) - \text{Lin}_{h,i}({y})|}{|x - y|}.
\]
Within a single region \( C_{h,i} \), since \( \text{Lin}_{h,i}(y) \) is linear, its Lipschitz constant is given by $L_{h,i} = |{a}_{h,i}|$,
where \( |{a}_{h,i}| \) is the norm of the gradient vector \( {a}_{h,i}\) of \( \text{Lin}_{h,i}({y}) \). The global Lipschitz constant of the piecewise linear approximation is
\[
L_{\text{pwl},h} = \max_{i=1,2,\dots,h} L_{h,i} = \max_{i=1,2,\dots,h} |{a}_{h,i}|.
\]
Each linear function \( \text{Lin}_{h,i}({y}) \) is constructed to approximate \( f(y) \) within cell \( C_{h,i} \). Typically
 \( \text{Lin}_{h,i}({y}) \) matches \( f(y) \) exactly at the vertices of \(  C_{h,i} \).
For example, in simplicial interpolation, the gradient is derived from the values of \( f({y}) \) at the vertices.
Since \( f({y}) \) is Lipschitz continuous with constant \( L \), the gradient \( {a}_{h,i} \) of \( \text{Lin}_{h,i}({y}) \) within \( C_{h,i} \) is derived from \( f(y) \). Specifically, \( |{a}_{h,i}| \) reflects the "local slope" of \( f(y) \) within \( C_{h,i}\). Since \( f(y) \) is Lipschitz continuous with constant \( L \), the maximum rate of change of \( f(y) \) is bounded by \( L \), that is
\[
|\nabla f(y)| \leq L, \quad \forall y \in \mathbb{R}^n.
\]
Because \( \text{Lin}_{h,i}({y}) \) is constructed to approximate \( f(y) \) and its gradient \( {a}_{h,i} \) is derived from \( f(y) \), we have
\[
|{a}_{h,i}| \leq L, \quad  i=1,2,\dots,h.
\]
The global Lipschitz constant of the piecewise linear approximation is the maximum gradient norm across all regions
\[
L_{\text{pwl},h} = \max_{i=1,2,\dots,h} |{a}_{h,i}|\leq L.
\]
To analyze the behavior of \( {a}_{h,i} \) (the gradient of the piecewise linear approximation) as \(h \to \infty \), we must consider how the partition of the domain \( \mathbb{R}^n \) evolves and how the approximation behaves as the number of regions \( C_{h,i} \) increases. 
Since the piecewise linear approximation is exact at the vertices of the simplex, \( f(y) = \text{Lin}_{h,i}({y}) \) at these points. This ensures that as \( h \to \infty \), the approximation becomes finer and better represents \( f(y) \).
Each \( {a}_{h,i} \) represents the slope of the linear approximation \( \text{Lin}_{h,i}({y}) \) within region \( C_{h,i} \), and it is determined by the values of \( f(y) \) at the vertices of \( C_{h,i} \). Specifically, \( {a}_{h,i} \) approximates the gradient of \( f(y) \) over \( C_{h,i} \). As \( h \to \infty \), the regions \( C_{h,i} \) become arbitrarily small, effectively shrinking to points. Within each region \( C_{h,i} \), the gradient \( {a}_{h,i} \) approximates the gradient, \( \nabla f(y) \), over that region, under assumption that $f(y)=f(y_1,y_2,\dots,y_n)$ is at least differentiable in $D$.  
\section[]{Existence and uniqueness of a solution for nonlinear system}
In this section, we formulate and prove the main results regarding the existence and uniqueness of the solution to the nonlinear system (\ref{nonlinear_system_n_order}) with a zero initial conditions. To achieve this, we employ a specially chosen weighted supremum norm, which ensures that the Laplace transform can be applied to the system. Subsequently, this approach allows us to formally define the transfer function of the nonlinear system as the limiting case of its piecewise linearization. This, in turn, facilitates the standard analysis of the system using transfer function algebra within a closed-loop control framework.
\begin{thm} Let $f\in C^2(\mathbb{R}^n\mapsto \mathbb{R})$ is globally Lipschitz in \((y, y', \dots, y^{(n-1)})\) with Lipschitz constant $L$ and $u\in C(\mathbb{R}\mapsto \mathbb{R})$. Then the initial value problem
\[
y^{(n)} + f(y, y', \dots, y^{(n-1)}) = u(t),
\]
with a zero initial conditions
\[
y(0) = y'(0) = \dots = y^{(n-1)}(0) = 0,
\]
has a unique solution $y(t)$ defined on $\mathbb{R}$ satisfying here the inequality
\[
|y(t)|\leq Ke^{2\tilde L|t|},\quad  K\geq0,\quad \tilde L=({n-1+L^2})^{1/2}.
\]
\end{thm}
\begin{pf}
The \(n\)-th order differential equation can be reformulated as a system of first-order ODEs (\ref{nonlinear_system_1_order}).
On the space \( C(\mathbb{R} \to \mathbb{R}^n) \) of continuous functions, define the weighted supremum norm
\[
\|y\|_{\tilde L} = \sup_{t \in \mathbb{R}} e^{-2\tilde L|t|}\,|y(t)|.
\]
The subspace of functions with finite norm under this definition includes all constant functions. Including all constant functions in this space ensures that simple initial guesses, such as \(y_0(t) = 0\), lie within the domain of the Picard operator, allowing the iterative process to start consistently within the function space.
The Picard operator
\( (P(y))(t) \) with zero initial conditions is defined as
   \[
   (P(y))(t) = \int_0^t \big(\tilde f(y(\tau))+\tilde u(\tau) \big)\, d\tau,
   \]
   where \( \tilde f \) is a Lipschitz function with Lipschitz constant \( \tilde L > 0 \). This means
   \[
   |\tilde f(x) - \tilde f(y)| \leq \tilde L |x - y|, \quad \forall x, y\in \mathbb{R}^n.
   \]
To prove that \( (P(y))(t)) \) is contractive, we show that the mapping \( P: y(t) \mapsto (P(y))(t) \) satisfies the contraction condition, that is, there exists \( 0 \leq k < 1 \) such that
\[
\|P(x) - P(y)\|_{\tilde L} \leq k \|x - y\|_{\tilde L}.
\]
By definition
\[
(P(x))(t) - (P(y))(t) = \int_0^t \big(\tilde f(x(\tau)) - \tilde f(y(\tau))\big) \, d\tau.
\]
The norm is given by  
\[
\|P(x) - P(y)\|_{\tilde L} 
= \sup_{t \in \mathbb{R}} e^{-2\tilde L |t|} \left| \int_0^t \big(\tilde f(x(\tau)) - \tilde f(y(\tau))\big) \, d\tau \right|.
\]
By the Lipschitz condition on \( \tilde f \),  
\[
|\tilde f(x(\tau)) - \tilde f(y(\tau))| \leq \tilde L |x(\tau) - y(\tau)|.
\]  
This implies 
\[
\left| \int_0^t \big(\tilde f(x(\tau)) - \tilde f(y(\tau))\big) \, d\tau \right| \leq \int_0^t \tilde L |x(\tau) - y(\tau)| \, d\tau.
\]
Using the definition of the norm \( \|\cdot\|_{\tilde L} \),  
\[
|x(\tau) - y(\tau)| \leq e^{2\tilde L |\tau|} \|x - y\|_{\tilde L}.
\]  
Substituting this bound into the integral yields
\[
\left| \int_0^t \big(\tilde f(x(\tau)) - \tilde f(y(\tau))\big) \, d\tau \right| \leq \int_0^t \tilde L e^{2\tilde L |\tau|} \|x - y\|_{\tilde L} \, d\tau.
\]
To express the estimate in terms of the weighted supremum norm, multiply both sides by the factor \( e^{-2\tilde L |t|} \) 
\[
e^{-2\tilde L |t|} \left| \int_0^t \big(\tilde f(x(\tau)) - \tilde f(y(\tau))\big) \, d\tau \right| 
\leq e^{-2\tilde L |t|} \int_0^t \tilde L e^{2\tilde L |\tau|} \|x - y\|_{\tilde L} \, d\tau.
\]
For \( \tau \in [0, t] \), we have \( |\tau| = \tau \). Hence,  
\[
e^{-2\tilde L |t|} e^{2\tilde L |\tau|} = e^{-2\tilde L (t - \tau)}.
\]
Substituting this expression, we obtain 
\[
e^{-2\tilde L |t|} \left| \int_0^t \big(\tilde f(x(\tau)) - \tilde f(y(\tau))\big) \, d\tau \right| 
\leq \int_0^t \tilde L e^{-2\tilde L (t - \tau)} \|x - y\|_{\tilde L} \, d\tau.
\]
The exponential term \( e^{-2\tilde L (t - \tau)} \) integrates to a finite value 
\[
\int_0^t e^{-2\tilde L (t - \tau)} \, d\tau = \frac{1 - e^{-2\tilde L t}}{2\tilde L} \leq \frac{1}{2\tilde L}.
\]
Therefore,
\[
e^{-2\tilde L |t|} \left| \int_0^t \big(\tilde f(x(\tau)) - \tilde f(y(\tau))\big) \, d\tau \right| 
\leq  \frac{1}{2} \|x - y\|_{\tilde L}.
\]
For \( t < 0 \), the argument proceeds analogously, with only minor adjustments due to the sign of \( t \). These involve appropriately handling the integration limits and the exponential weighting term \( e^{-2\tilde L |t|} \), but the underlying reasoning remains identical.

Taking the supremum over \( t \in \mathbb{R} \), we obtain
\[
\|P(x)- P(y)\|_{\tilde L} \leq \frac{1}{2} \|x - y\|_{\tilde L}.
\]
Since \( \frac{1}{2} < 1 \), the operator \( P \) is a contraction. 

To apply the Banach Fixed-Point Theorem, we first need to verify that the Picard operator \( P(y) \) maps the chosen function space into itself. Specifically, we must show that if \( y(t) \) belongs to the space of continuous functions endowed with the weighted supremum norm \( \|y\|_{\tilde L} = \sup_{t \in \mathbb{R}} e^{-2\tilde L |t|} |y(t)| < \infty \), then \( P(y)(t) \) also lies in the same space.

Define the space \( \mathcal{C}_{\tilde L} \) as
   \[
   \mathcal{C}_{\tilde L} = \{ y \in C(\mathbb{R} \mapsto \mathbb{R}^n) \; : \; \|y\|_{\tilde L} < \infty \}.
   \]
A function \( y(t) \in \mathcal{C}_{\tilde{L}} \) satisfies
\[
\|y\|_{\tilde{L}} = \sup_{t \in \mathbb{R}} e^{-2\tilde{L} |t|} |y(t)| < \infty.
\]
This implies that, for all \( t \in \mathbb{R} \),
\[
|y(t)| \leq K e^{2\tilde{L}|t|},
\]
where \( K = \|y\|_{\tilde{L}} \) is a finite constant.

Recall that the Picard operator is defined by
\[
(P(y))(t) = \int_0^t \big(\tilde{f}(y(\tau)) + \tilde u(\tau)\big) \, d\tau,
\]
where \( \tilde{f} \) satisfies the Lipschitz condition with constant \( \tilde{L} > 0 \), and \( \tilde u(t) \in \mathcal{C}_{\tilde L} \) is a continuous input function.  
We aim to show that
\[
\text{if } y \in \mathcal{C}_{\tilde L}, \text{ then } P(y) \in \mathcal{C}_{\tilde L}.
\]

Since \( y \in \mathcal{C}_{\tilde L} \), it follows that \( \|\tilde{f}(y)\|_{\tilde L} < \infty \) because \( \tilde{f} \) is Lipschitz continuous, and thus
\[
\|\tilde{f}(y)\|_{\tilde L} \leq |\tilde{f}(0)| + \tilde{L} \|y\|_{\tilde L}.
\]
Moreover, \( \tilde u(t) \in \mathcal{C}_{\tilde L} \) is continuous and satisfies \( \|\tilde u\|_{\tilde L} < \infty \).  
Hence, the integrand \( \tilde{f}(y(\tau)) + \tilde u(\tau) \) is well defined and integrable for all \( t \in \mathbb{R} \).
To complete the proof, we need to verify that
\[
\|P(y)\|_{\tilde L} = \sup_{t \in \mathbb{R}} e^{-2\tilde L |t|} 
\left|\int_0^t \big(\tilde{f}(y(\tau)) + \tilde u(\tau)\big) \, d\tau \right| < \infty.
\]
Using the triangle inequality, we obtain
\[
\left|\int_0^t \big(\tilde{f}(y(\tau)) + \tilde u(\tau)\big) \, d\tau \right| 
\leq \int_0^t |\tilde{f}(y(\tau))| \, d\tau + \int_0^t |\tilde u(\tau)| \, d\tau.
\]
Applying the weighted supremum norm to both sides yields
\[
e^{-2\tilde{L}|t|} \left|\int_0^t \big(\tilde{f}(y(\tau)) + \tilde u(\tau)\big) \, d\tau \right| 
\leq e^{-2\tilde{L}|t|} \int_0^t e^{2\tilde{L}|\tau|} \|\tilde{f}(y)\|_{\tilde L} \, d\tau 
+ e^{-2\tilde{L}|t|} \int_0^t e^{2\tilde{L}|\tau|} \|\tilde u\|_{\tilde L} \, d\tau.
\]
Consequently,
\[
e^{-2\tilde{L}|t|} \left|\int_0^t \big(\tilde{f}(y(\tau)) + \tilde u(\tau)\big) \, d\tau \right|
\leq \|\tilde{f}(y)\|_{\tilde L} \cdot \frac{1}{2\tilde{L}} (1 - e^{-2\tilde{L}|t|}) +\|\tilde u\|_{\tilde L} \cdot \frac{1}{2\tilde{L}} (1 - e^{-2\tilde{L}|t|}).
\]
Taking the supremum over \( t \in \mathbb{R} \), we conclude that the norm is finite,
\[
\|P(y)\|_{\tilde L} \leq \frac{\|\tilde{f}(y)\|_{\tilde L} +\|\tilde u\|_{\tilde L} }{2\tilde{L}} < \infty.
\]
\end{pf} 
\begin{rem}
The Laplace transform of \(y(t)\) defined as
\[
\mathcal{L}_-[y(t)] = \int_{0^-}^\infty y(t)e^{-st} \, dt\triangleq Y(s),
\]
where \(s \in \mathbb{C}\) with \(\Re(s) = \sigma > 0\) is well-defined for all functions \(y(t) \in \mathcal{C}_{\tilde{L}}\) provided \(\Re(s) = \sigma > 2\tilde{L}\). This condition ensures that the exponential decay induced by \(e^{-\sigma t}\) dominates any growth in \(y(t)\), making the integral converge.

In the definition, \( t = 0^- \) refers to the limit where \( t \) approaches \( 0 \) from the left-hand side, capturing the pre-initial values of the function \( y(t) \). 
This framework allows the analysis of functions that may exhibit discontinuities or impulsive behavior, such as those involving the Dirac delta \( \delta(t) \) and its derivative. The Dirac delta \( \delta(t) \) is a generalized function (distribution) defined by
\[
\int_{-\infty}^{\infty} \delta(t)\, dt = 1, \qquad \delta(t) = 0 \quad \text{for} \quad t \neq 0.
\]

This definition of the Laplace transform implies the time-derivative rule
\[
\mathcal{L}_-[y'(t)] = sY(s) - y(0^-),
\]
whereby initial conditions existing before $t = 0 $ are brought into the analysis. This formulation ensures that the Laplace transform correctly incorporates the effects of any discontinuities at \(t = 0\).

For higher derivatives, the rule generalizes as
\[
\mathcal{L}_-[y^{(n)}(t)] = s^nY(s) - s^{n-1}y(0^-) - \cdots - y^{(n-1)}(0^-).
\]
Zero initial conditions  
\[
y(0) = y'(0) = \dots = y^{(n-1)}(0) = 0,
\]  
when using the Laplace transform, will be understood in the sense of  
\[
y(0^-) = y'(0^-) = \dots = y^{(n-1)}(0^-) = 0.
\]  

The initial-value theorem states that 
\begin{equation}\label{initialValueThm}
\lim_{s \to \infty \cdot 1} sY(s) = y(0^+),
\end{equation}
where \(s \to \infty \cdot 1\) indicates that the limit is taken along the positive real axis. This theorem is particularly useful for validating solutions of differential equations and determining the immediate response of a system to initial conditions. An excellent and inspiring discussion on the topic of the Laplace transform can be found in the paper~\cite{Lundberg2007}.

In the context of solving linear differential equations using the Laplace transform, the pre-initial values \(y(0^-)\), \(y'(0^-), \dots\)  plays a critical role in converting the problem into the Laplace domain. By incorporating both pre-initial and post-initial values, this approach ensures that discontinuities in \(y(t)\) or its derivatives are properly handled and  impulsive inputs such as \(\delta(t)\) or \(\delta'(t)\) are accurately incorporated into the solution.
\end{rem}
\section[]{Convergence of the solution of the piecewise linear approximation to the solution of the original nonlinear system}
Consider the infinite sequence of piecewise linear approximation of the original nonlinear system (\ref{nonlinear_system_n_order})
\begin{equation}\label{piecewise_linear_system_n_order}
y_h^{(n)}(t) +  \text{Lin}_{h}({\bar y_h(t)}) = {u}_{h}(t),\quad t\geq0
\end{equation}
with zero initial conditions
\[
y_h(0) = y_h'(0) = \dots = y_h^{(n-1)}(0) = 0,
\]
where $\bar y_h(t)\triangleq\big(y_h(t),y'_h(t),\dots y^{(n-1)}_h(t) \big)$, and $\{h \}$ is an infinite and increasing sequence of natural numbers,
given that \( u_{h}(t) \to u(t) \) uniformly on every interval \( [0,T]\subset \mathbb{R}\), $\text{Lin}_{h}({y})$ is globally Lipschitz with the same constant $L$ as $f(y)$, and ${u}_{h}, u\in \mathcal{C}_{\tilde L}$. Analogously as in (\ref{nonlinear_system_1_order}), for equivalent first-order system we use notation
\begin{equation}\label{eq_LinSys}
y'_h=\tilde{\text{Lin}_{h}}({y_h}) + \tilde{u}_{h}(t)
\end{equation}
As has been proven before $\tilde{\text{Lin}_{h}}(y)\in\mathcal{C}_{\tilde L}$, \( \text{Lin}_{h}({y}) \) converges to \( f(y) \) uniformly  on any compact domain $D\subset\mathbb{R}^n$ if $f$ is twice continuously differentiable. 
Thus, for any compact set \( D \subset \mathbb{R}^n \) and any \(\epsilon_f > 0\), there exists \( h_0 \) such that for all \( h \geq h_0 \)
\[
\sup_{y \in D} |\text{Lin}_{h}({y}) - f(y)| < \epsilon_f.
\]
To study the relationship between \( y_h(t) \) and \( y(t) \), define the difference
\[
z_h(t) = y_h(t) - y(t).
\]
Subtract the equations (\ref{eq_LinSys}) and (\ref{nonlinear_system_1_order}) we have
\[
\frac{dz_n}{dt} = \left[\tilde{\text{Lin}_{h}}({y_h}) - \tilde f(y)\right] + \left[\tilde u_h(t) - \tilde u(t)\right].
\]
Using the triangle inequality
\[
|\tilde{\text{Lin}_{h}}({y_h}) - \tilde f(y)| 
\leq |\tilde{\text{Lin}_{h}}({y_h}) - \tilde{\text{Lin}_{h}}({y})| + |\tilde{\text{Lin}_{h}}({y}) - \tilde f(y)|.
\]
The term \( |\tilde{\text{Lin}_{h}}({y_h}) - \tilde{\text{Lin}_{h}}({y})| \) depends on the Lipschitz continuity of \( \tilde{\text{Lin}_{h}} \)
   \[
   |\tilde{\text{Lin}_{h}}({y_h}) - \tilde{\text{Lin}_{h}}({y})| \leq \tilde L_h |y_h - y| = \tilde L_h |z_h|,
   \]
where \( \tilde L_h \) is the Lipschitz constant for \( \tilde{\text{Lin}_{h}}(y) \). For large \( h \), \( \tilde L_h \) can be bounded by a common Lipschitz constant \( \tilde L \) for \( \tilde f(y) \), due to uniform convergence. The term \( |\tilde{\text{Lin}_{h}}({y}) - \tilde f(y)| \) depends on the uniform convergence of \( \tilde{\text{Lin}_{h}}(y) \) to \( \tilde f(y) \). On any compact \( D \), this term can be made arbitrarily small by choosing \( h \) sufficiently large. Thus
\[
\left|\frac{dz_h(t)}{dt}\right| \leq \tilde L |z_h(t)| + (\epsilon_f + \epsilon_u),
\]
where \( \epsilon_f > 0 \) and $\epsilon_u >0$ accounts for the uniform convergence error \( |\tilde{\text{Lin}_{h}}({y}) - \tilde f(y)| \) on $D$ and \( |\tilde{u_{h}}({t}) - \tilde u(t)| \) on $[0,T]$, respectively.  Using the inequality
\[
\frac{d}{dt}|z_h(t)|\leq\left|\frac{d}{dt}z_h(t)\right|
\]
and Gronwall inequality (Theorem~\ref{Gronwall}), we can bound \( |z_h(t)| \)
\[
|z_h(t)| \leq (\epsilon_f + \epsilon_u) \frac{e^{\tilde Lt} - 1}{\tilde L}.
\]

For any fixed \( T\) and given \( (\epsilon_f + \epsilon_u) > 0 \), there exists \( h_0 \) such that \( h \geq h_0\) ensures \( |z_h(t)| \) is arbitrarily small. This implies uniform convergence of \( y_h(t) \to y(t) \) on the finite interval $[0,T]$. The following theorem summarizes the preceding results.
\begin{thm}\label{thm:convergence}
Consider the infinite sequence of piecewise linear approximations of the original nonlinear system (\ref{nonlinear_system_n_order}):  
\[
y_h^{(n)}(t) + \text{Lin}_{h}(\bar{y}_h(t)) = u_h(t), \quad t \geq 0,
\]  
where 
\(
\bar{y}_h(t) \triangleq \big(y_h(t), y_h'(t), \dots, y_h^{(n-1)}(t)\big),
\)
with zero initial conditions, and \( \{h\} \) is an infinite, increasing sequence of natural numbers. Assume that  
 \( u_h(t) \to u(t) \) uniformly on any finite interval \( [0, T] \subset \mathbb{R} \) and  \( u_h(t), u(t) \in \mathcal{C}_{\tilde{L}} \).  

Then, for any compact set \( D \subset \mathbb{R}^n \) and any finite \( T > 0 \), the solution \( y_h(t) \) of the piecewise linearized system converges uniformly to the solution \( y(t) \) of the original nonlinear system as \( h \to \infty \):  
\[
\sup_{t \in [0, T]} |y_h(t) - y(t)| \to 0.
\]
\end{thm}
\section[]{Definition of transfer function for nonlinear system}
Building on the convergence result established in the previous section, 
we now turn our attention to the piecewise linear approximation of the original nonlinear system.
For the system represented by (\ref{piecewise_linear_system_n_order}),
\[
y_h^{(n)}(t) + \text{Lin}_h(\bar{y}_h(t)) = u_h(t), \quad t \ge 0,
\]
we proceed by transforming the equation into the frequency domain. 
Taking the Laplace transform of both sides, we obtain
\[
\mathcal{L}_-\!\big[y_h^{(n)}(t)\big] + \mathcal{L}_-\!\big[\text{Lin}_h(\bar{y}_h(t))\big] = \mathcal{L}_-\!\big[u_h(t)\big].
\]
Assuming zero initial conditions \( y_h(0) = y_h'(0) = \dots = y_h^{(n-1)}(0) = 0 \) and using standard properties of the Laplace transform, we have
\[
\mathcal{L}_-\!\big[y_h^{(n)}(t)\big] = s^n Y_h(s),
\]
where \( Y_h(s) \) denotes the Laplace transform of \( y_h(t) \). Substituting this expression into the previous equation gives
\[
s^n Y_h(s) + \mathcal{L}_-\!\big[\text{Lin}_h(\bar{y}_h(t))\big] = U_h(s),
\]
where \( U_h(s) = \mathcal{L}_-\![u_h(t)] \).

If \( \text{Lin}_h(\bar{y}_h) = a_h \bar{y}_h + b_h \), where \( a_h \) is the cumulative gradient of the linear function in compact domain $D$, its Laplace transform becomes
\[
\mathcal{L}_- \big[ \text{Lin}_h(\bar{y}_h) \big] = a_h \cdot (1,s,\dots,s^{n-1}) Y_h(s) + \mathcal{L}_-[b_h].
\]
Substituting, we have
\[
s^n Y_h(s) + a_h \cdot (1,s,\dots,s^{n-1}) Y_h(s) + \frac{b_h}{s} = U_h(s).
\]
Reorganizing the equation,
\[
\big(s^n  + a_h \cdot (1,s,\dots,s^{n-1}) \big) Y_h(s) = U_h(s) - \frac{b_h}{s}.
\]
Thus, the solution becomes
\[
Y_h(s) = \frac{1}{\big(s^n  + a_h \cdot (1,s,\dots,s^{n-1})\big)} \left(U_h(s) - \frac{b_h}{s}\right)
\]
and the transfer function \( G_h(s) \) of the system is 
\[
G_h(s) = \frac{1}{s^n  + a_h \cdot (1,s,\dots,s^{n-1})},
\]
where \(\cdot\) stands for the standard dot product. However, the input now includes a term due to \( b_h \)
\[
Y_h(s) = G_h(s) U_h(s) - G_h(s) \frac{b_h}{s}.
\]
The transfer function describes the relationship between the input \( u_h(t) \) and the output \( y_h(t) \) in the Laplace domain. For piecewise linear systems, \( G_h(s) \) may vary across different regions \( C_{h,i} \), depending on the coefficients \( a_{h,i} \) and \( b_{h,i} \).

As \( h \to \infty \), the piecewise linear approximations \( \mathrm{Lin}_h({y}) \) converge uniformly to the nonlinear function \( f(y) \) on compact domains. Similarly, \( u_h(t) \to u(t) \) uniformly on any finite interval \( [0, T] \). According to Theorem~\ref{thm:convergence}, this implies that the solution \( y_h(t) \) of the linearized system converges uniformly to the solution \( y(t) \) of the original nonlinear system \(y^{(n)}(t) + f(\bar{y}(t)) = u(t)\) as \( h \to \infty \).

For this nonlinear system, we cannot directly define a transfer function due to the nonlinearity of \( f(y) \). However, in the limit of piecewise linearized regime, the transfer function takes the form
\begin{equation}\label{transfer:nonlinear}
G(s) = \frac{1}{s^n  + \nabla f(y) \cdot (1,s,\dots,s^{n-1})}.
\end{equation}
Here, \( \nabla f(y) \) represents the gradient of \( f(y)=f(y_1,y_2,\dots,y_n) \) at the operating point $y\in D\subset\mathbb{R}^n$.

Consider a feedback control system, as depicted schematically in Figure~\ref{fig1_control}, 
\[
\begin{aligned}
    Y_h(s) &= G_h(s) U_h(s) - G_h(s) \frac{b_h}{s}, \\
    U_h(s) &= C_h(s) E_h(s), \\
    E_h(s) &= R(s) - Y_h(s).
\end{aligned}
\]
The closed-loop system equation is given as
\[
\big(1 + C_h(s) G_h(s)\big) Y_h(s) = C_h(s)G_h(s) R(s) - G_h(s) \frac{b_h}{s},
\]
where

\( C_h(s) = K_{p,h} + \frac{K_{i,h}}{s} + K_{d,h} s \) is the transfer function of the PID controller,

\( G_h(s) \) is the plant's piecewise linearization transfer function,

\( R(s) = \frac{1}{s} \) is the Laplace transform of the Heaviside unit step function \( \eta(t) \), defined as 
\[
\eta(t) = 
\begin{cases}
	0, & t < 0, \\
	1, & t \ge 0,
\end{cases} 
\]

\( b_h \) acts like a source of persistent input that does not decay (if \( b_h \neq 0 \)).

Substituting \( G_h(s), C_h(s) \) and $R(s)$ we obtain
\[
\left(s^n  + a_h \cdot (1,s,\dots,s^{n-1})+K_{p,h} + \frac{K_{i,h}}{s} + K_{d,h} s\right)Y_h(s)
= \left(K_{p,h} + \frac{K_{i,h}}{s} + K_{d,h} s\right)\frac1s-\frac{b_h}{s}
\]
Now multiply through by \( s \) to clear the denominators
\[
\big(s^{n+1} +  a_h \cdot (s, s^2, \dots, s^{n}) +  K_{p,h}s + K_{i,h} +  K_{d,h}s^2\big) Y_h(s)
=K_{p,h} + \frac{K_{i,h}}{s} +  K_{d,h}s - b_h.
\]
In the time-domain form of the differential equation under zero initial conditions, obtained by applying the inverse Laplace transform, each term on the left-hand side corresponds to derivatives of \(y_h(t)\):

\begin{itemize}
	\item \(s^{n+1} Y_h(s)\) corresponds to the \((n+1)\)-th derivative \(y_h^{(n+1)}(t)\),
	\item \(a_h \cdot (s, s^2, \dots, s^n) Y_h(s)\) corresponds to terms up to the \(n\)-th derivative, namely\newline \(a_h \cdot (y_h'(t), y_h''(t), \dots, y_h^{(n)}(t))\),
	\item \(K_{p,h} s Y_h(s)\) corresponds to \(K_{p,h} y_h'(t)\),
	\item \(K_{d,h} s^2 Y_h(s)\) corresponds to \(K_{d,h} y_h''(t)\),
	\item \(K_{i,h} Y_h(s)\) corresponds directly to \(K_{i,h} y_h(t)\).
\end{itemize}

Hence, the left-hand side in the time domain becomes
\[
y_h^{(n+1)}(t) + a_h \cdot (y_h'(t), y_h''(t), \dots, y_h^{(n)}(t)) + K_{p,h} y_h'(t) + K_{i,h} y_h(t) + K_{d,h} y_h''(t).
\]

The right-hand side in the Laplace domain is
\[
K_{p,h} + \frac{K_{i,h}}{s} + K_{d,h}s - b_h.
\]

Applying the inverse Laplace transform term by term yields:

\begin{itemize}
	\item \(K_{p,h}\) corresponds to \(K_{p,h} \delta(t)\), a Dirac delta impulse,
	\item \(\frac{K_{i,h}}{s}\) corresponds to \(K_{i,h} \eta(t)\), with \(\eta(t)\) being the Heaviside unit step function,
	\item \(K_{d,h}s\) corresponds to \(K_{d,h} \delta'(t)\), where \(\delta'(t)\) is the Dirac doublet (the first derivative of the Dirac delta in the sense of distributions),
	\item \(b_h\) corresponds to \(b_h \delta(t)\).
\end{itemize}

It follows from the general rule for the Laplace transform of the \(n\)-th derivative of the Dirac delta, \(\delta^{(n)}(t)\), understood in the distributional sense, that
\[
\mathcal{L}_-[\delta^{(n)}(t)] = s^n,
\]
see \cite{papoulis1962} and \cite{Lundberg2007}.  
Therefore, the right-hand side in the time domain becomes
\[
(K_{p,h} - b_h)\delta(t) + K_{i,h} \eta(t) + K_{d,h} \delta'(t).
\]
Now, combining all terms, the time-domain differential equation reads
\begin{equation}\label{model:general}
\begin{aligned}
y_h^{(n+1)}(t) + a_h \cdot (y_h'(t), y_h''(t), \dots, y_h^{(n)}(t))
+ K_{p,h} y_h'(t) + K_{i,h} y_h(t) + K_{d,h} y_h''(t) \\
= \left(K_{p,h}- b_h  \right)\delta(t) + K_{i,h} \eta(t) + K_{d,h} \delta'(t),
\end{aligned}
\end{equation}
where the left-hand side represents the dynamics of the linearized system under the PID controller, and the right-hand side includes impulsive contributions \( \delta(t) \) and \( \delta'(t) \) arising from the proportional and derivative actions, as well as the persistent input \(K_{i,h} \eta(t)\).

Consider now the special case \( n = 1 \). Then, equation~\eqref{model:general} reduces to
\begin{equation}\label{eq:first_order}
	\begin{aligned}
		(1 + K_{d,h}) y_h''(t) + (a_h + K_{p,h}) y_h'(t) + K_{i,h} y_h(t) 
		&= \left(K_{p,h}- b_h\right) \delta(t) + K_{i,h} \eta(t) + K_{d,h} \delta'(t),
	\end{aligned}
\end{equation}
with zero initial conditions
\begin{equation}\label{eq:first_orderIC}
	y_h(0^-) = y'_h(0^-) = 0.
\end{equation}

Under the assumptions that \( f \in C^2(\mathbb{R}\mapsto\mathbb{R}) \) and has a globally bounded first derivative, the solutions \( y_h(t) \) of the piecewise linear approximations converge to the solution \( y(t) \) of the original nonlinear system
\[
y'(t) + f(y(t)) = u(t)
\]
as \( h \to \infty \).

\begin{figure}[ht]
	\begin{center}
		\includegraphics[width=\textwidth]{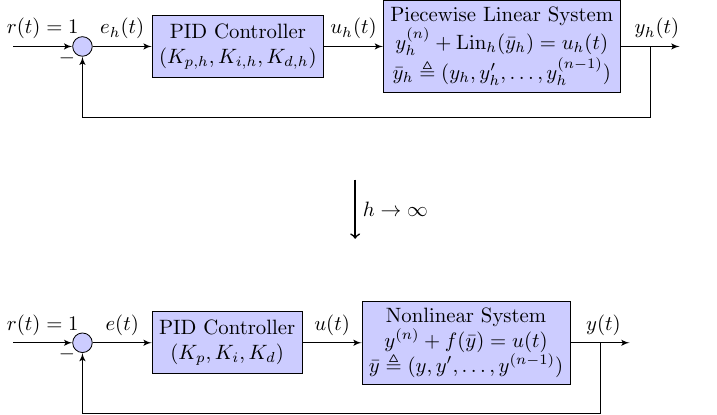}    
		\caption{Convergence of the piecewise linear approximation to the original nonlinear control system.}  
		\label{fig1_conv}                                  
	\end{center}                                 
\end{figure}

If the piecewise linear approximations converge uniformly, the interchange of the limit and the Laplace transform is justified, ensuring that the frequency-domain analysis of the approximated models accurately reflects the behavior of the original system, as illustrated in Figure~\ref{fig1_conv}. 

\section[]{Numerical simulations}
In this section, we focus on numerical simulations to illustrate the effectiveness of the proposed method based on piecewise linearization of nonlinear systems. Both simulations consider first-order systems, using the model given in (\ref{eq:first_order}) and (\ref{eq:first_orderIC}). This formulation explicitly incorporates the zero initial state into the output response and, at the same time, mitigates issues arising from jump and impulse discontinuities on the right-hand side of the differential equation by employing approximations of the impulsive functions \( \delta(t) \) and \( \delta'(t) \). This approach is consistent with the discussion in Remark~\ref{rem:nonzero}.

To approximate the Dirac delta and its derivative (the Dirac doublet), we employ a narrow Gaussian function
\[
\delta(t) \approx \frac{1}{\sigma\sqrt{2\pi}} e^{-\frac{t^2}{2\sigma^2}},
\]
where \( \sigma \) is a small positive value that controls the width of the Gaussian peak. As \( \sigma \) becomes smaller, the Gaussian function becomes more sharply peaked at \( t = 0 \), and its integral over all time approaches $1$, mimicking the behavior of the Dirac delta. To approximate the Dirac doublet \( \delta'(t) \), you would differentiate the Gaussian approximation of \( \delta(t) \). The derivative of the Gaussian function with respect to \( t \) is
\[
\delta'(t) \approx -\frac{t}{\sigma^2} \cdot \frac{1}{\sigma \sqrt{2\pi}} e^{-\frac{t^2}{2\sigma^2}},
\]
where \( \sigma \) is again a small positive value that controls the width of the Gaussian peak. 

Dirac doublet \( \delta'(t) \) is a distribution that behaves like a pulse with a singularity at \( t = 0 \) and is widely used in systems where impulses and their derivatives are involved, such as in control systems, signal processing, or electromagnetics.

\begin{figure}[ht]
    \centering
    \begin{minipage}{0.45\textwidth}
        \centering
        \includegraphics[width=\textwidth]{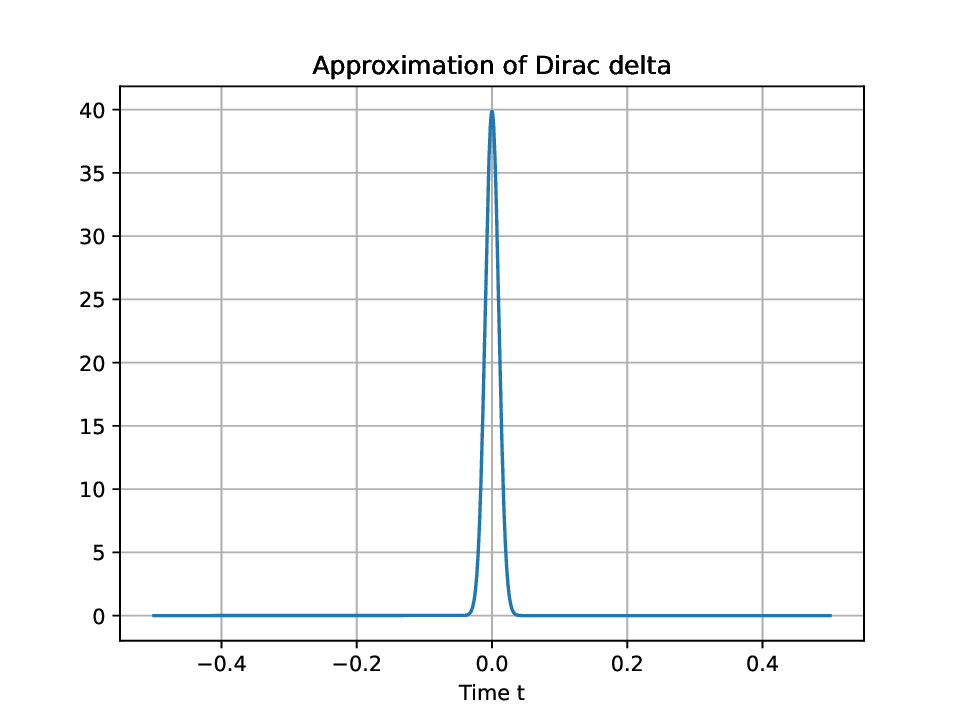}
        \caption{Gaussian approximation of the Dirac delta ($\delta(t)$) with standard deviation $\sigma = 0.01$. }
        \label{fig:fig1D}
    \end{minipage}\hfill
    \begin{minipage}{0.45\textwidth}
        \centering
        \includegraphics[width=\textwidth]{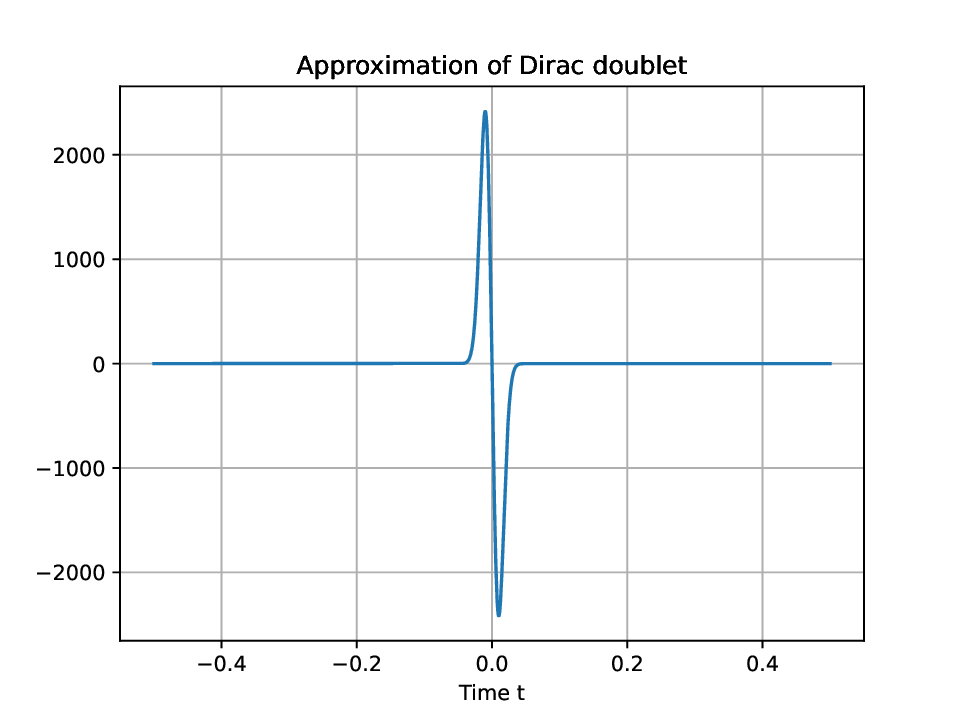}
        \caption{Gaussian approximation of the Dirac doublet ($\delta'(t)$) with standard deviation $\sigma = 0.01$.}
        \label{fig:fig2DD}
    \end{minipage}
\end{figure}

In all simulations, the Dirac delta \( \delta(t) \) and its derivative \( \delta'(t) \) are approximated using Gaussian functions with \( \sigma = 0.01 \). For clarity, both approximations are illustrated in Figures~\ref{fig:fig1D} and~\ref{fig:fig2DD}. Since the Gaussian approximation does not rely on a limiting process as \( h \to \infty \), it provides a fixed, finite-width representation suitable for computational purposes. This approach ensures numerical stability and feasibility in simulations and represents a commonly used method when Dirac delta functions and their derivatives are required in continuous-time numerical models.

To optimize the PID controller parameters \( (K_{p,h}, K_{i,h}, K_{d,h}) \) using the PSO method, we define the cost function as a combination of two widely used criteria: the ITAE (Integral of Time-weighted Absolute Error) and  the ISO (Integral of Squared Overshoot). The cost function 
\( J(K_{p,h}, K_{i,h}, K_{d,h}) \) is expressed as a weighted sum of these two criteria.  The ITAE criterion is defined as
\[
\text{ITAE}(K_{p,h}, K_{i,h}, K_{d,h}) = \int_0^\infty t \, |e_h(t)| \, dt,
\]
where \( e_h(t) = r(t) - y_h(t) \) represents the error between the reference signal \( r(t) \) and the system output \( y_h(t) \). The ITAE penalizes both the magnitude and the duration of the error, favoring responses that reduce the error more quickly.

The ISO criterion penalizes overshoot in the system's response, limiting the peak amplitude above the desired reference during transients. It is defined as
\[
\text{ISO}(K_{p,h}, K_{i,h}, K_{d,h}) = \int_0^\infty \left( \max(0, y_h(t) - r(t)) \right)^2 dt,
\]
where the term \( \max(0, y_h(t) - r(t)) \) captures positive deviations above the reference, and the square ensures that larger overshoots are penalized more heavily. The combined cost function is then
\[
J(K_{p,h}, K_{i,h}, K_{d,h}) 
= \lambda_{\text{ITAE}} \, \text{ITAE}(K_{p,h}, K_{i,h}, K_{d,h})
+ \lambda_{\text{ISO}} \, \text{ISO}(K_{p,h}, K_{i,h}, K_{d,h}),
\]
where the weights \( \lambda_{\text{ITAE}} \) and \( \lambda_{\text{ISO}} \) control the relative importance of each criterion during optimization.

In practice, it is often unnecessary to integrate over the infinite time domain. Instead, the integrals are evaluated over a finite time interval \([0, T]\) that captures the relevant system behavior, such as the settling period. Thus, the cost function used in the PSO optimization becomes
\begin{equation}\label{ITAE_ISO_criterion}
J(K_{p,h}, K_{i,h}, K_{d,h}) 
= \int_0^T \left[ t \, |r(t)-y_h(t)| + \alpha \left( \max(0, y_h(t) - r(t)) \right)^2 \right] dt,
\end{equation}
where \( \alpha > 0 \) balances overshoot reduction against error minimization, and \( r(t) = \eta(t) \) is the Heaviside unit step function.  
This finite-time formulation is commonly employed in numerical simulations and optimization algorithms such as PSO to efficiently determine optimal PID controller parameters.
\subsection{Example 1: PID control for a linear first-order system}
The first example, based on a linear system, is included to verify the correctness of the proposed optimization procedure and to demonstrate consistency with classical control results. This provides a baseline validation before extending the same approach to nonlinear systems.

Let us consider a first-order dynamic system defined by
\begin{equation}\label{linear:control}
\frac{dy(t)}{dt} + a\,y(t) = b\,u(t),
\end{equation}
where \(y(t)\) denotes the system output, \(u(t)\) is the control input, \(a > 0\) represents the system decay rate (the reciprocal of the time constant), and \(b > 0\) is the system gain. The goal is to design a PID controller that ensures the system output \(y(t)\) tracks the desired reference signal \(r(t) = \eta(t)\). 

The corresponding transfer function of the plant in the Laplace domain is
\[
G(s) = \frac{Y(s)}{U(s)} = \frac{b}{s + a}.
\]
For instance, with parameters \(a = 2\) and \(b = 1\), the transfer function becomes
\[
G(s) = \frac{1}{s + 2},
\]
which represents a stable first-order system with a time constant 
\(\tau = \frac{1}{a} = 0.5\).

PID controller in the Laplace domain is expressed as
\[
C(s) = K_p + \frac{K_i}{s} + K_d s,
\]
and the corresponding closed-loop transfer function is
\[
T(s) = \frac{C(s)G(s)}{1 + C(s)G(s)}.
\]
The controller parameters \((K_p, K_i, K_d)\) are selected to achieve the desired closed-loop performance specifications. For a first-order plant, a straightforward design procedure based on closed-loop tuning rules \cite{Astrom1995PID} can be employed. Using the standard tuning relations
\[
K_p = \frac{1.2}{b\tau}, \qquad K_i = \frac{2}{\tau}, \qquad K_d = \frac{0.5\tau}{b},
\]
and substituting \(b = 1\) and \(\tau = 0.5\), we obtain
\[
K_p = 2.4, \qquad K_i = 4.0, \qquad K_d = 0.25.
\]
Hence, the resulting PID controller is
\[
C(s) = 2.4 + \frac{4.0}{s} + 0.25s.
\]

\begin{figure}[ht]
    \centering
    \begin{minipage}{0.49\textwidth}
        \centering
        \includegraphics[width=\textwidth]{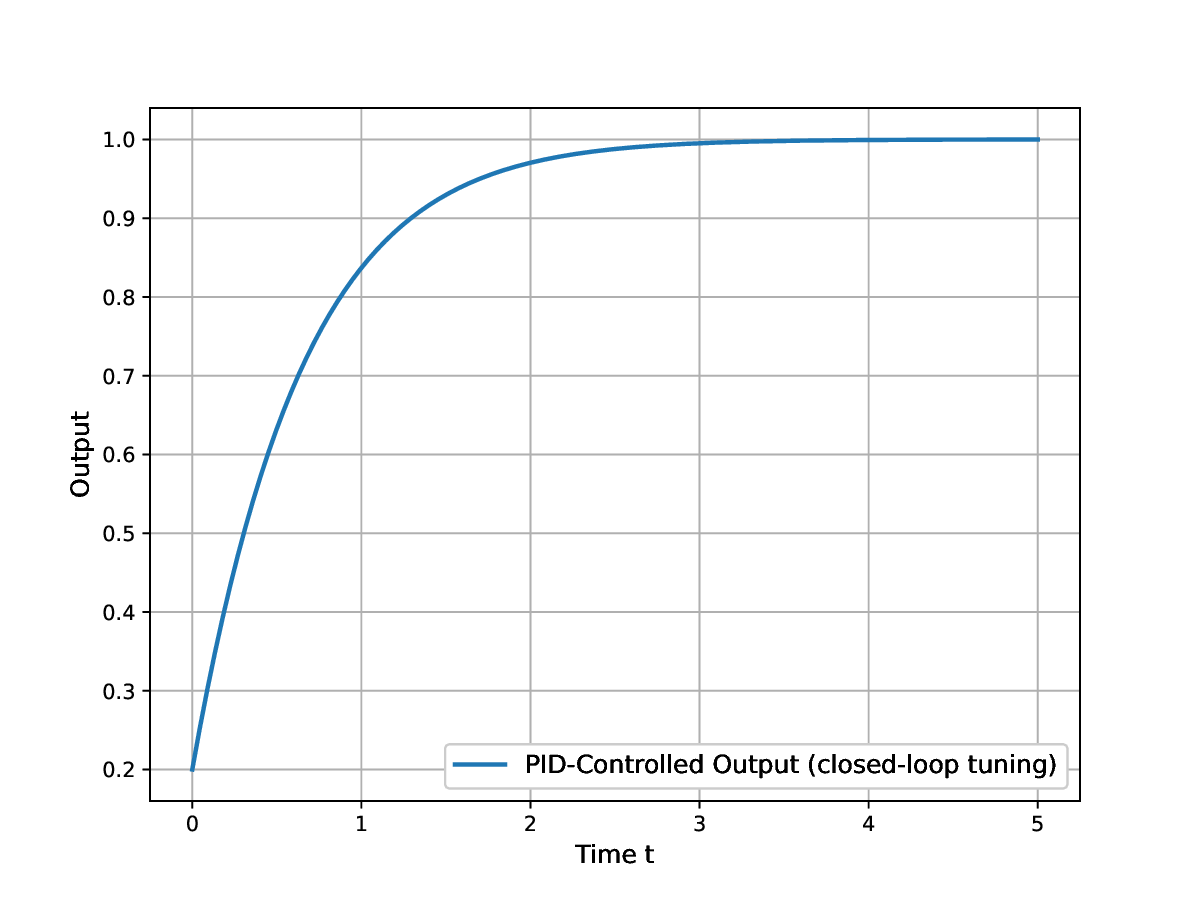}
        \caption{$C(s) = 2.4 + \frac{4.0}{s} + 0.25s$. 
PID-controlled output obtained using the closed-loop tuning method. 
The plot illustrates the time evolution of the system response under the computed PID parameters, highlighting the transient behavior and the tracking of the reference signal. This simulation provides a benchmark for comparing the performance of alternative controller tuning approaches.}

        \label{fig:fig1ex1}
    \end{minipage}\hfill
    \begin{minipage}{0.49\textwidth}
        \centering
        \includegraphics[width=\textwidth]{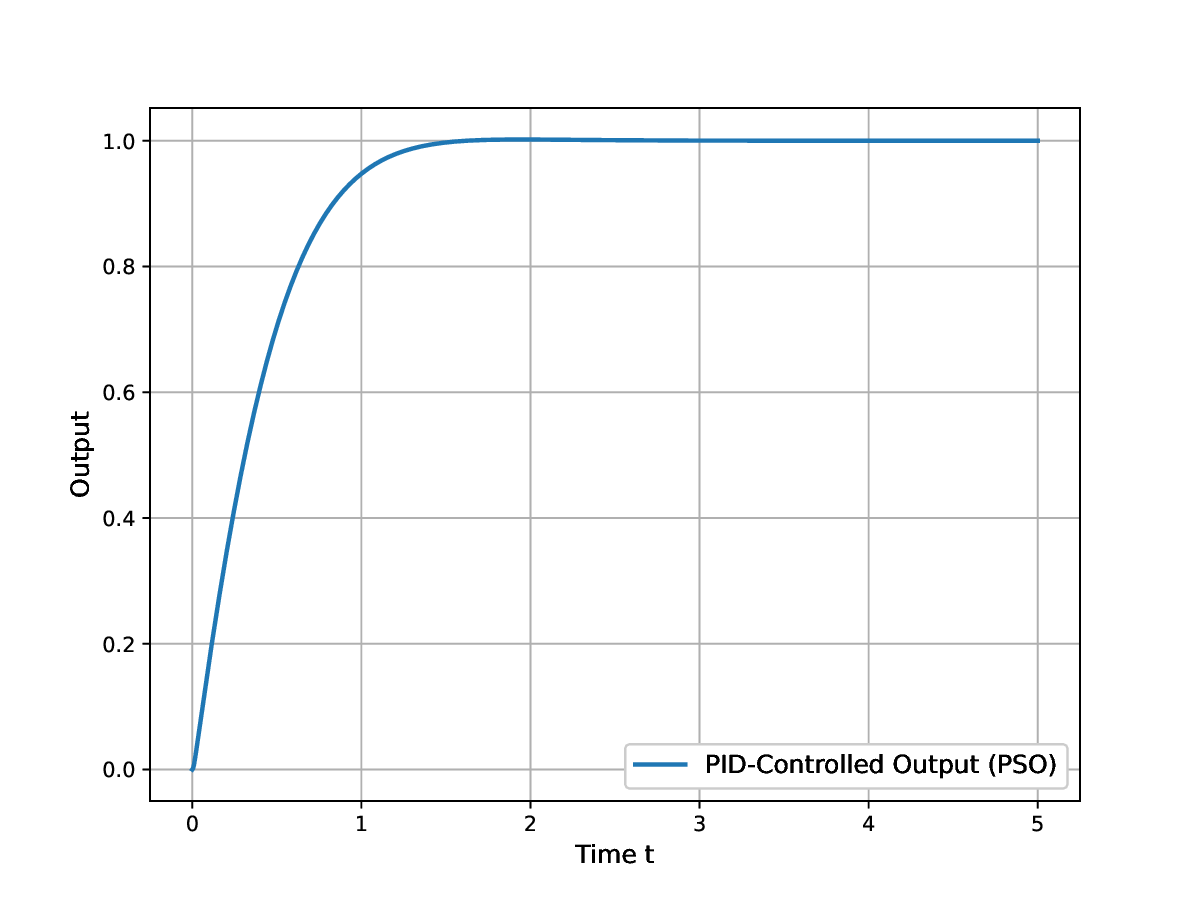}
        \caption{$C(s) = 3.72 + \frac{10.0}{s} + 0.0s$. 
PID-controlled output (PSO-based tuning), using the combined ITAE+ISO criterion (\ref{ITAE_ISO_criterion}) with $\alpha = 2000$, obtained after 5 iterations with a swarm size of 30, for 
$(K_{p,h}, K_{i,h}, K_{d,h}) \in [0, 10]^3$.       
}
\label{fig:fig1ex1PID}
    \end{minipage}
\end{figure}

The control flow comparison for the PID controller parameters, determined using the closed-loop tuning method and the PSO algorithm is presented in Figures~\ref{fig:fig1ex1} and~\ref{fig:fig1ex1PID}.
\begin{rem}
The system under consideration is linear; hence, $\text{Lin}_{h}({y})=f({y})=ay$, for ${y} \in \mathbb{R}$ and any chosen partition parameter $h$.
\end{rem}
\begin{rem}\label{rem:nonzero}
For the computed controller parameters \( K_p = 2.4 \), \( K_i = 4.0 \), and \( K_d = 0.25 \), 
the corresponding unit-step response of the closed-loop system---given by (\ref{eq:first_order}), 
and, in the case of the linear model (\ref{linear:control}), detailed in Appendix~\ref{app:IJandIS}---is described by the solution of the second-order differential equation
\[
1.25y'' + 4.4y' + 4y = 0.25\delta'(t) + 2.4\delta(t) + 4\eta(t),
\]
which can be derived using the Laplace transform under the initial conditions 
\( y(0^-) = 0 \) and \( y'(0^-) = 0 \).

However, certain subtleties in the treatment of distributional terms, particularly those involving the Dirac delta and its derivative, may lead to apparent inconsistencies between the obtained solution, the prescribed initial conditions, or the original differential equation itself. A careful interpretation in the sense of distributions resolves these discrepancies and ensures consistency of the solution.

These issues are typically related to how the Dirac delta and its derivative are handled in the Laplace transform process. When \( \delta(t) \) and \( \delta'(t) \) appears on the right-hand side, it modifies the initial conditions of the ODE. Specifically, the term \( 0.25\delta'(t) \) introduces a discontinuity in \( y'(t) \), so the initial derivative \( y'(0^+) \) is no longer \( 0 \), even if \( y'(0^-) = 0 \). Similarly, \( \delta(t) \) modifies \( y(0^+) \), making it nonzero. Thus, the initial conditions after \( t = 0 \) (that is, \( y(0^+) \) and \( y'(0^+) \)) must be carefully derived by integrating the ODE across \( t = 0 \). The Dirac delta and its derivative have the following properties: 
   \[
   \int_{0^-}^{0^+} \delta(t) \, dt = 1.
   \]
This is because the Dirac delta integrates to $1$ over any interval containing \(t = 0\). For the derivative of the Dirac delta
   \[
   \int_{0^-}^{0^+} \delta'(t) \, dt = \delta(t) \Big|_{0^-}^{0^+}=0
   \]
   because \(\delta(t)\) is defined to be zero outside \(t = 0\).
   
Now analyze left-hand side across $t=0$. The derivative \(y''\) contributes discontinuities in \(y'(t)\),
\[
     \int_{0^-}^{0^+} y''(t) \, dt = y'(0^+) - y'(0^-),
\]
analogously, the derivative \(y'(t)\) contributes discontinuities in \(y(t)\),
\[
     \int_{0^-}^{0^+} y'(t) \, dt = y(0^+) - y(0^-),
\]
and \(\int y(t)\) is continuous, so
     \[
     \int_{0^-}^{0^+} y(t) \, dt = 0.
     \]
Substituting and equating both sides, we obtain
\begin{equation}\label{eq:postInit}
1.25y'(0^+) + 4.4y(0^+) = 2.4.
\end{equation}
From the unit step response of this system using the transfer function, 
\[
Y(s) =\frac{T(s)}{s} = \frac{K_d s^2 + K_p s + K_i}{(1 + K_d)s^3 + (2 + K_p)s^2 + K_i s}
= \frac{0.25s^2 + 2.4s + 4.0}{1.25s^3 + 4.4s^2 + 4.0s}.
\]
Using the initial-value theorem (\ref{initialValueThm}), the post-initial value is obtained as 
\( y(0^+) = \lim_{s \to \infty} sY(s) = 0.2 \).
Thus, from (\ref{eq:postInit}), \( y'(0^+) = 1.216 \).

Taking the inverse Laplace transform of \( Y(s) \), the time-domain representation of the unit step response is
\[
y(t) = \mathcal{L}^{-1}\{Y(s)\} 
= 1 - \frac{4}{5} e^{-\frac{44t}{25}} \cos\!\left(\frac{8t}{25}\right) 
- \frac{3}{5} e^{-\frac{44t}{25}} \sin\!\left(\frac{8t}{25}\right),
\]
but as can be easily verified, this function does not satisfy the differential equation nor the initial conditions.  
Indeed, the function \( y(t) \) satisfies the homogeneous form of the equation,
\[
1.25y'' + 4.4y' + 4y = 4,
\]
which corresponds to the regular (non-distributional) part of the system response after the impulsive effects have vanished.  
The Dirac delta \( \delta(t) \) and its derivative \( \delta'(t) \) appearing in the original equation 
\[
1.25y'' + 4.4y' + 4y = 0.25\delta'(t) + 2.4\delta(t) + 4\eta(t)
\]
act only at \( t = 0 \), producing instantaneous changes in \( y \) and \( y' \).  
After these impulsive contributions are accounted for, the evolution for \( t > 0 \) is governed by the smooth equation above.  
The correct solution that satisfies both the equation and the initial conditions is therefore 
\[
y_{\text{true}}(t) = y(t)\eta(t),
\]
where \( \eta(t) \) denotes the Heaviside unit step function, ensuring that the response starts at \( t = 0 \) with the appropriate jump conditions induced by the impulses.

This fact can be rigorously verified in the sense of distributions by applying the standard rules for a smooth function \( g(t) = y(t) \), see~\cite{Lundberg2007},
\[
\begin{aligned}
	\eta'(t) &= \delta(t), \quad \eta''(t) = \delta'(t), \\
	\delta(t)g(t) &= \delta(t)g(0), \\
	\delta'(t)g(t) &= \delta'(t)g(0) - \delta(t)g'(0).
\end{aligned}
\]
Accordingly, the derivatives of \( y_{\text{true}}(t) \) can be written as
\[
\begin{aligned}
	y_{\text{true}}'(t) &= y'(t)\eta(t) + y(0)\delta(t), \\
	y_{\text{true}}''(t) &= y''(t)\eta(t) + y'(0)\delta(t) + y(0)\delta'(t).
\end{aligned}
\] 
Substituting these expressions into the differential equation gives
\[
1.25y_{\text{true}}'' + 4.4y_{\text{true}}' + 4y_{\text{true}}
\]
\[
= [1.25y'' + 4.4y' + 4y]\eta(t) 
+ [1.25y'(0) + 4.4y(0)]\delta(t) + 1.25y(0)\delta'(t).
\]
Matching the coefficients of \( \eta(t) \), \( \delta(t) \), and \( \delta'(t) \) yields 
\( y(0^+) = 0.2 \) and \( y'(0^+) = 1.216 \), 
which confirms that \( y_{\text{true}}(t) \) satisfies both the differential equation 
\[
1.25y'' + 4.4y' + 4y = 0.25\delta'(t) + 2.4\delta(t) + 4\eta(t),
\]
and the prescribed initial conditions \( y(0^-)=y'(0^-)=0 \). 
\end{rem}
\begin{rem}\label{initialJumpandSlope}
This method represents a generally applicable approach for analyzing \(n\)-th order nonlinear systems by utilizing the model (\ref{model:general}) and the principles of infinitesimal calculus as outlined in Remark~\ref{rem:nonzero}. 

For the special case \(n = 1\), the solution yields the following expressions (the derivation of these relations is provided in Appendix~\ref{app:IJandIS}).

The initial jump at \(t = 0^+\) is given by
\[
y(0^+) = \frac{b K_d}{1 + b K_d},
\]
indicating that the immediate response of the system depends on the derivative 
term, which can introduce a discontinuity in the output.

The initial slope of the response at \(t = 0^+\) is
\[
y'(0^+) =  \frac{b K_p (1 + b K_d) - b K_d (a + b K_p)}{(1 + b K_d)^2},
\]
representing the instantaneous rate of change in the system output immediately after the input is applied.

Thus, the exact solution of equation (\ref{eq:first_order}), applied to Example 1 with initial conditions \( y(0^-) = 0 \), \( y'(0^-) = 0 \), and PID settings as specified in Figure 6 (\( K_p = 3.72 \), \( K_i = 10 \), \( K_d = 0 \)), yields a zero initial jump at \( t = 0^+ \) due to the absence of a derivative component (\( K_d = 0 \)). Furthermore, the initial slope of the response is \( y'(0^+) = 3.72 \), corresponding to a rise angle of \( \phi = \arctan(y'(0^+)) \approx 75^\circ \).

These results offer valuable insight into the transient behavior of the system, capturing the immediate influence of the control parameters on the output 
response following the application of the input.
\end{rem}
\subsection{Example 2: PID control for a nonlinear first-order system}\label{example_nonlinear}
In this example, we demonstrate the practical application of the proposed methodology to a nonlinear first-order system. The objective is not only to assess closed-loop performance, but also to illustrate a typical scenario in which the piecewise linear approximation provides a natural and effective framework for PID controller design.

In many practical control problems, the exact analytical form of the nonlinear function \( f(\cdot) \) governing the system dynamics is not fully known. Instead, the system behavior is often characterized through experimental measurements, steady-state operating points, or tabulated input--output data. Under such conditions, constructing a globally valid nonlinear model is difficult, whereas local linear approximations over a compact operating domain are readily available. The piecewise linear model employed in this example can therefore be interpreted as an interpolation of known system data, rather than an artificial linearization of a fully known nonlinear equation.

The PID controller parameters are optimized using the piecewise linear approximation, which enables the use of standard performance criteria and efficient numerical optimization. Due to the uniform convergence of the piecewise linear model and the corresponding system responses, the resulting controller is subsequently applied without modification to the original nonlinear system. The simulation results confirm that this strategy yields closed-loop behavior comparable to that obtained from direct nonlinear tuning, while retaining analytical transparency and computational efficiency.

Consider the nonlinear first-order system  
\begin{equation*}
y'(t) + f(y(t)) = u(t), \quad t \geq 0,
\end{equation*}
where the nonlinear function \( f(y) \) is defined as  
\begin{equation*}
f(y) = 0.5y + \ln(1 + y^2).
\end{equation*}
The function \( f(y) \) belongs to the class \( C^2(\mathbb{R} \to \mathbb{R}) \), which ensures it is twice continuously differentiable over the entire real line. Moreover, \( f(y) \) satisfies the following properties:
\begin{itemize}
    \item \( f(y) \to \infty \) as \( |y| \to \infty \), indicating that the function exhibits unbounded growth for large values of \( y \).
    \item Both the first and second derivatives of \( f \) are globally bounded on \( \mathbb{R} \). Specifically,
    \[
    |f'(y)| \leq 1.5, \quad |f''(y)| \leq 2, \quad \forall y \in \mathbb{R}.
    \]
\end{itemize}
Let the interval of linearization \( D \subset \mathbb{R} \) be chosen as \( D = [-3, 3] \).  
If, during the simulation, the output \( y_h \) of the closed-loop control system approaches the boundary of the domain \( D \), the interval will be expanded as needed.

\begin{figure}[htbp]
		\centering
		
		\begin{minipage}{0.45\textwidth}
			\centering
			\includegraphics[width=\linewidth,height=5cm,keepaspectratio]{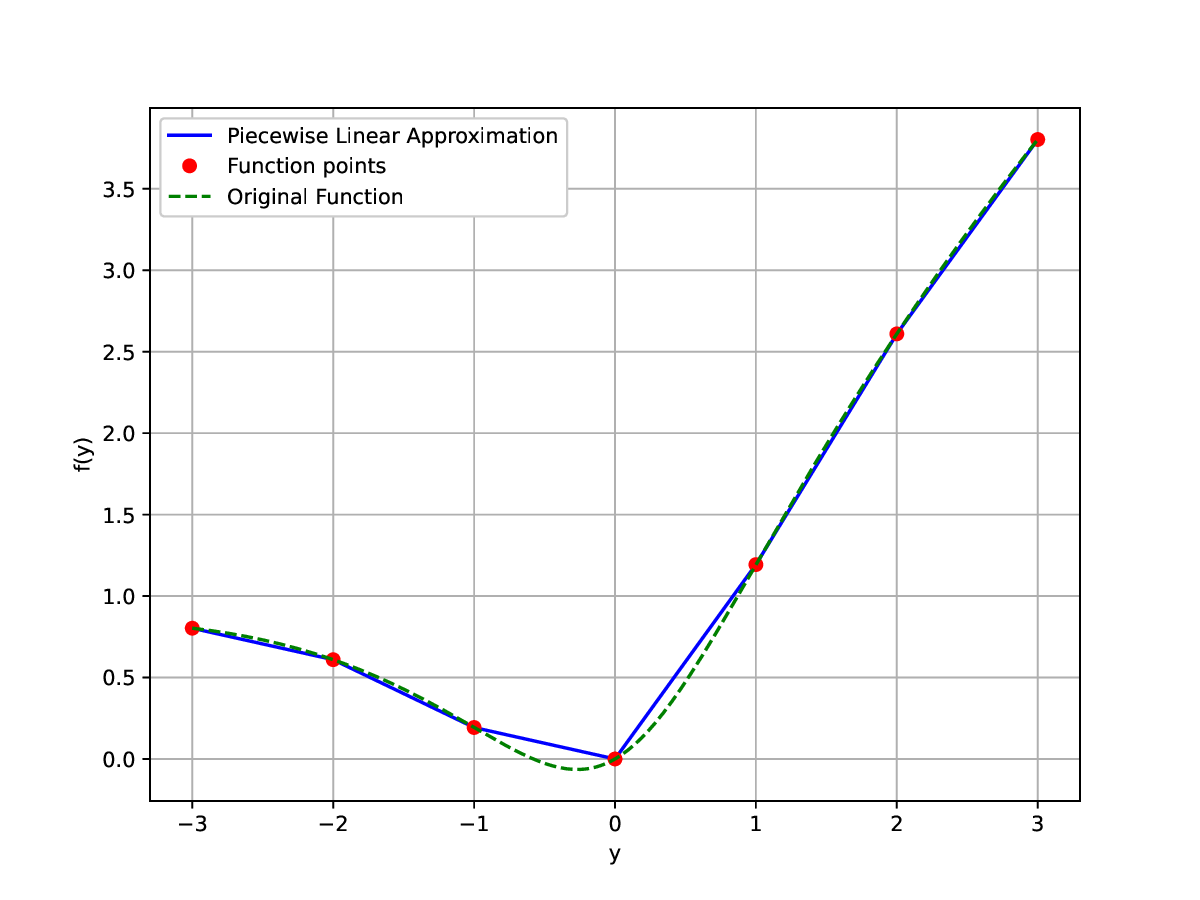}
			\caption{The nonlinear function $f(y) = 0.5y + \ln(1 + y^2)$ and its piecewise linear approximation $\text{Lin}_{h=6}(y)$ used in the construction of the piecewise linear model.}  
			\label{fig_ex2}        
		\end{minipage}
		
		\vspace{1em} 
		
		\begin{minipage}{0.45\textwidth}
			\centering
			\includegraphics[width=\linewidth,height=4.7cm,keepaspectratio]{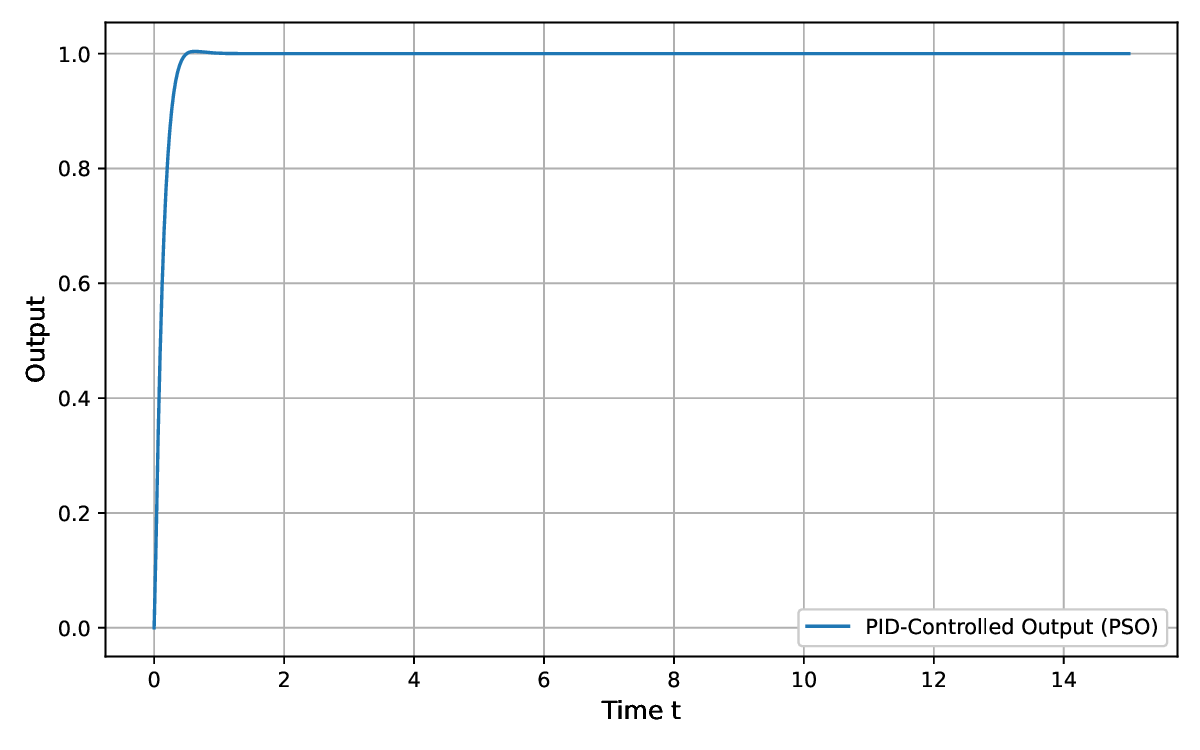}
			\caption{$C(s) = 13.5927 + \frac{65.1469}{s} + 0.0s$. 
				PID-controlled output (PSO-based tuning), using the combined ITAE+ISO criterion (\ref{ITAE_ISO_criterion}) with $\alpha = 5$, obtained after 10 iterations with a swarm size of 30, for 
				$(K_{p,6}, K_{i,6}, K_{d,6}) \in [0, 100]^3$.}  
			\label{fig_Ex2} 
		\end{minipage}
		\hspace{1em} 
		\begin{minipage}{0.45\textwidth}
			\centering
			\includegraphics[width=\linewidth,height=5cm,keepaspectratio]{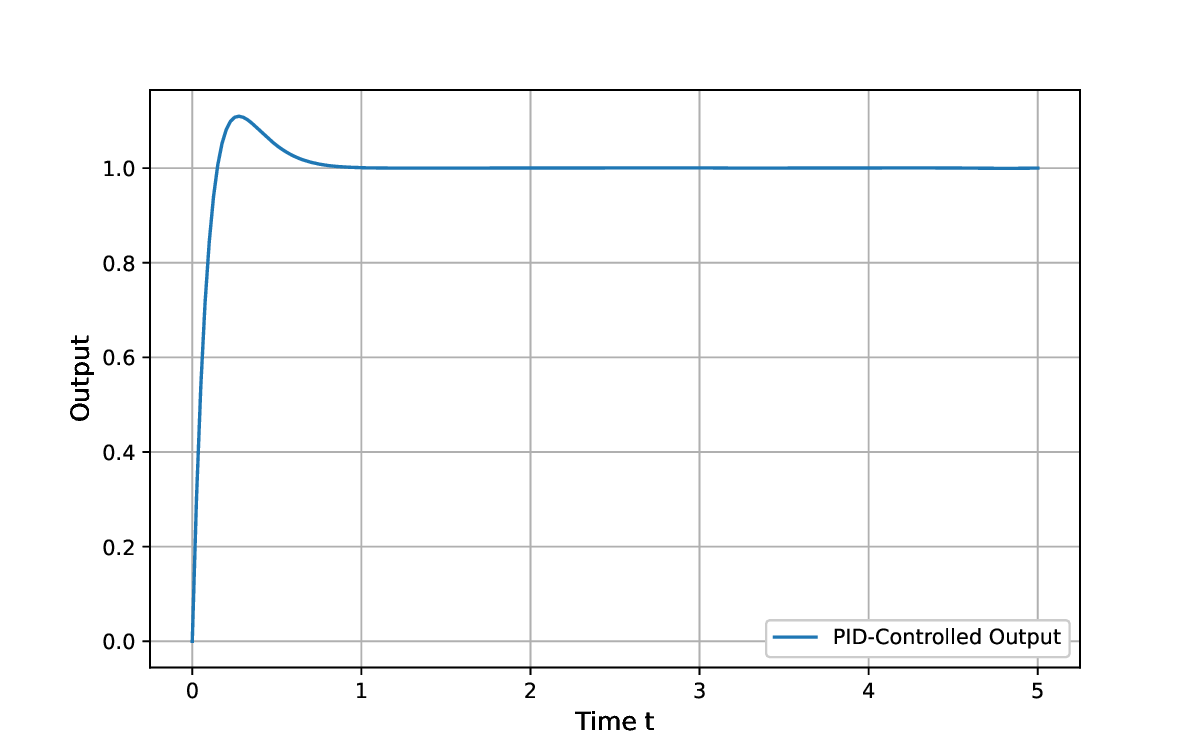}
			\caption{Simulation of the nonlinear first-order system with $f(y) = 0.5y + \ln(1 + y^2)$ controlled by a PID regulator according to the scheme in Figure~\ref{fig1_conv}, with PID parameters \(K_{p} = 13.5927\), \(K_{i} = 65.1469\), \(K_{d} = 0.0\), obtained via PSO-based tuning for the piecewise linear approximation $\text{Lin}_{6}(y)$ of the original nonlinear system.}  
			\label{fig_Ex2nonlinear}
		\end{minipage}
		
		\vspace{1em} 
		
		\begin{minipage}{0.45\textwidth}
			\centering
			\includegraphics[width=\linewidth,height=5cm,keepaspectratio]{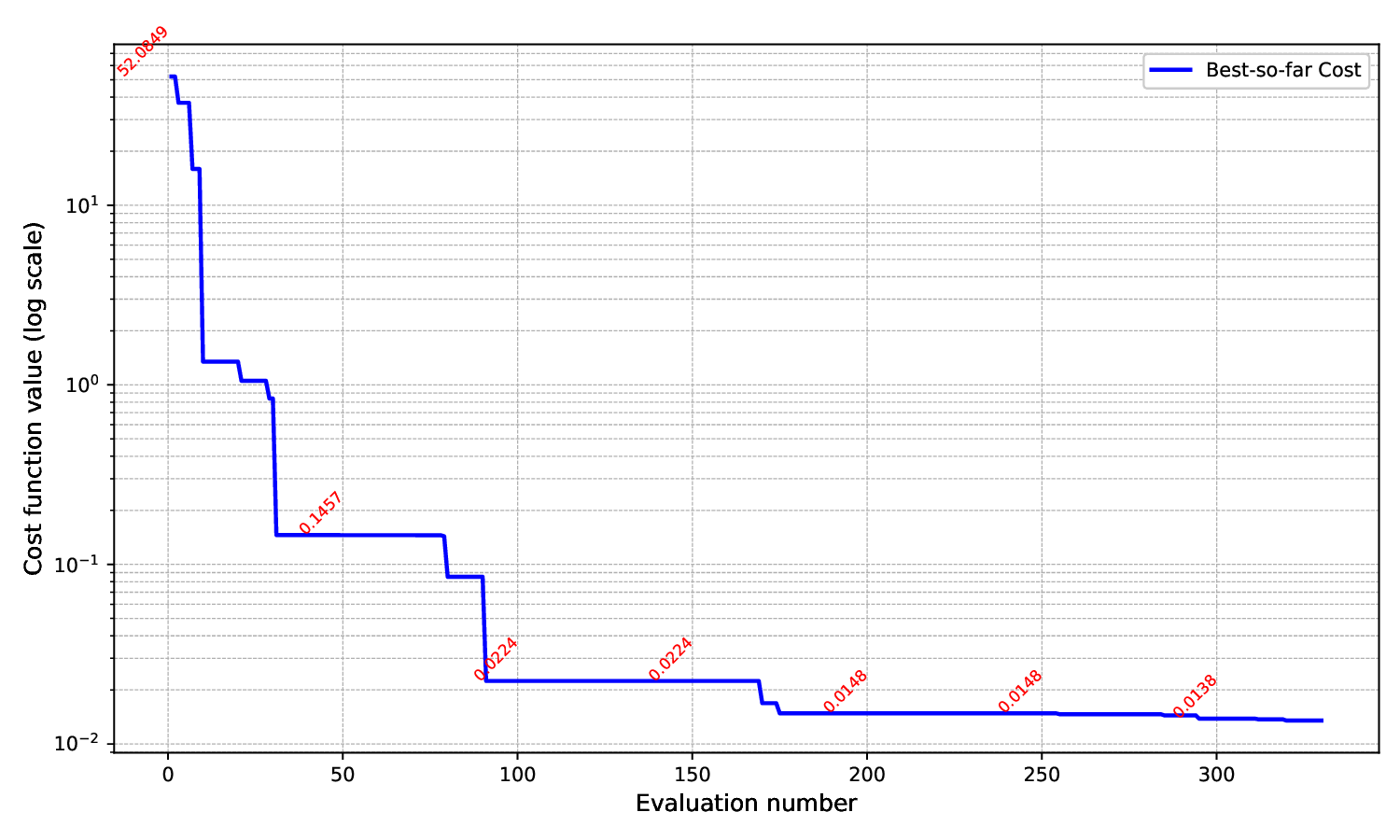}  
			\caption{Evolution of the PSO cost function over successive evaluations, showing the running minimum ("best-so-far") and selected annotated values to illustrate PID tuning progress.}
			\label{fig_Ex2nonlinear_cost}
		\end{minipage}
		
	\end{figure}

Assume \( h = 6 \), meaning that \( D \) is divided into six subintervals.  
The piecewise linear approximation of the nonlinear function
\[
f(y) = 0.5y + \ln(1 + y^2)
\]
is defined on subdomains \( C_{6,i} \) as summarized in Table~\ref{tab:Lin6},  
see also Figure~\ref{fig_ex2}. The complete piecewise linear approximation is then expressed as  
\[
\text{Lin}_{6}(y) \triangleq \{ \text{Lin}_{6,i}(y) \mid i = 1, 2, \dots, 6 \}.
\]
\begin{table}[ht]
\centering
\caption{Piecewise linear approximation $\text{Lin}_{6,i}(y)$ of $f(y)$ over $D=[-3,3]$.}
\label{tab:Lin6}
\begin{tabular}{ccl}
\toprule
Segment $i$ & Subinterval $C_{6,i}$ & Linear function $\text{Lin}_{6,i}(y)$ \\ 
\midrule
1 & $[-3, -2]$ & $\text{Lin}_{6,1}(y) = -0.19(y + 3) + 0.80$ \\
2 & $[-2, -1]$ & $\text{Lin}_{6,2}(y) = -0.42(y + 2) + 0.61$ \\
3 & $[-1, \phantom{-}0]$ & $\text{Lin}_{6,3}(y) = -0.19(y + 1) + 0.19$ \\
4 & $[\phantom{-}0, \phantom{-}1]$ & $\text{Lin}_{6,4}(y) = 1.19(y - 0) + 0.00$ \\
5 & $[\phantom{-}1, \phantom{-}2]$ & $\text{Lin}_{6,5}(y) = 1.42(y - 1) + 1.19$ \\
6 & $[\phantom{-}2, \phantom{-}3]$ & $\text{Lin}_{6,6}(y) = 1.19(y - 2) + 2.61$ \\
\bottomrule
\end{tabular}
\end{table}

The simulation results for the piecewise linear approximation and the original nonlinear system are presented in Figures~\ref{fig_Ex2} and~\ref{fig_Ex2nonlinear}, respectively. Figure~\ref{fig_Ex2} illustrates the PID-controlled output of the piecewise linear system with the controller parameters 
\[
C(s) = 13.5927 + \frac{65.1469}{s} + 0.0s,
\] 
which were obtained via PSO-based tuning using the combined ITAE+ISO criterion (\ref{ITAE_ISO_criterion}) with \(\alpha = 5\). Figure~\ref{fig_Ex2nonlinear} depicts the response of the original nonlinear first-order system controlled by a PID regulator with the same parameters derived from the piecewise linear approximation \(\text{Lin}_{6}(y)\). 

Figure~\ref{fig_Ex2nonlinear_cost} presents the evolution of the cost function over successive evaluations, providing insight into the convergence behavior of the PSO algorithm, the minimum achieved cost, and the total number of evaluations. To clearly illustrate the tuning process, only the running minimum of the cost function (that is, the "best-so-far" value) is plotted, effectively demonstrating the progressive improvement of the PID parameters throughout the iterations. Selected cost values are annotated at specific intervals to convey quantitative information while avoiding visual clutter.

Both simulations exhibit similar behavior, with small overshoot and a rapid approach to the reference value. The slightly larger overshoot observed in the nonlinear system can be attributed to the finite resolution of the piecewise linear approximation (\(h = 6\)), which introduces small discrepancies in the system dynamics compared to the true nonlinear model. Increasing the number of linear segments \(h\) would improve the approximation accuracy and further reduce these discrepancies, making the response of the nonlinear system even closer to that of the piecewise linear model. Despite these differences, the overall transient performance remains satisfactory, demonstrating the effectiveness of the PID parameters obtained via PSO tuning on the approximate model when applied to the original nonlinear system.


The linearization interval \(D\) is selected systematically based on the expected operating range of the system output. This ensures that the linear approximation remains valid across all practically relevant operating conditions. 

The number of subintervals \(h\) is determined analytically from the curvature of the nonlinear function \(f(y)\) using the standard interpolation error bound
\[
|f(y) - \mathrm{Lin}_h(y)| \le \tfrac{1}{2}\|\nabla^2 f\| \, \mathrm{diam}^2_{\max},
\]
where \(\mathrm{diam}_{\max}\) denotes the maximum width of an individual subinterval. This relation provides a direct link between the maximum subinterval width and the desired approximation accuracy \(\varepsilon\):
\[
\mathrm{diam}_{\max} \le \sqrt{\frac{2\varepsilon}{\|\nabla^2 f\|}}, \qquad
h = \left\lceil \frac{|D|}{\mathrm{diam}_{\max}} \right\rceil.
\]
If, during the simulation, the output of the closed-loop control system approaches the boundary of the domain \(D\), the interval \(D\) is expanded accordingly to ensure that the linear approximation remains valid throughout the entire operating region.  

This procedure provides a reproducible and computationally efficient way to define and maintain the parameter ranges without additional trial-and-error tuning or repeated simulations.

Expanding \( D \) dynamically ensures that the simulation remains valid for a wide range of system responses. However, it also increases computational complexity. Choosing an appropriate initial domain \( D \) is important to minimize unnecessary expansions. Dividing \( D \) into 6 equal subintervals simplifies the linearization process while maintaining a good approximation of the original nonlinear function \( f(y) \). Increasing \( h \) (more subdivisions) would improve the accuracy of the approximation but would also add computational overhead.
Since \( f(y) \) is being replaced by its piecewise linear approximation \( \text{Lin}_h(y) \), the system retains a globally Lipschitz structure.

The optimal controller parameters are obtained using PSO, where the objective function evaluates the closed-loop tracking performance of the nonlinear (piecewise linear) system. The optimization routine iteratively simulates the system response and minimizes the integral performance index. The general procedure of the proposed PSO-based tuning for first-order systems is outlined in Algorithm~~\ref{alg:pso_tuning}.

\begin{algorithm}[H]
\caption{PSO-based optimal tuning of controller parameters for a piecewise linear system}
\begin{algorithmic}[1]
    \State \textbf{Given:} piecewise linear model 
    \[
    {z}' = f_{\mathrm{PL}}(t, z, p), \quad 
    z = [y, {y}']^\top,
    \]
    with initial conditions \( z(0) = 0 \) 
    \hfill $\leftarrow$ (\ref{eq:first_order}), (\ref{eq:first_orderIC})
    
    \State Define the performance index:
    \[
    J(p) = \int_0^T t\,|1 - y(t, p)|\,dt 
    + \alpha \int_0^T (\max(0,\, y(t,p) - 1))^2\,dt
    \]
    
    \State Set search bounds for controller parameters 
    \( p = [p_1, p_2, p_3] = [K_p, K_i, K_d] \)
    
    \State Initialize a swarm of \( N \) particles 
    \( \{p^{(i)}\}_{i=1}^{N} \) with random velocities
    
    \For{each iteration \( k = 1, 2, \ldots, k_{\max} \)}
        \For{each particle \( p^{(i)} \)}
            \State Simulate the piecewise linear model 
            \({z}' = f_{\mathrm{PL}}(t, z, p^{(i)})\) over \(t \in [0, T]\)
            \State Compute the objective value \(J(p^{(i)})\)
            \State Update the personal and global best positions
        \EndFor
    \EndFor
    
    \State \textbf{Return:} optimal controller parameters 
    \( p^\ast = \arg\min_p J(p) \)
    
    \State Validate the closed-loop response \( y(t, p^\ast) \) on the piecewise linear model to confirm that the desired transient and steady-state performance is achieved.
    
    \State Apply the computed optimal controller parameters \( p^\ast \) to the original nonlinear system to evaluate the control performance in the full nonlinear context.
\end{algorithmic}
\label{alg:pso_tuning}
\end{algorithm}
\noindent
\textbf{PSO setup:} 
The optimization procedure was implemented using the \texttt{pyswarm} library in Python, employing the standard particle swarm optimization (PSO) algorithm with default internal parameters. 
The swarm consisted of $N = 30$ particles and the optimization was performed over 10 iterations. 
Each particle searched within the predefined parameter bounds 
\[
p = [p_1, p_2, p_3] \in [0,\,100]^3.
\]
The performance index was evaluated over a finite time horizon of $T = 15$, corresponding to the complete transient response of the closed-loop system. 
Default inertia weight and acceleration coefficients of the \texttt{pyswarm} implementation were used, as they provided stable convergence for the considered problem. 
The stopping criterion was based on the maximum number of iterations, which was sufficient since the best objective value stabilized well before this limit. 
Using the default configuration of the \texttt{pyswarm} package ensures full reproducibility of the results. 

As we might have noticed, for first-order systems without time delay, the derivative component is not required to achieve the desired setpoint -- a PI controller is typically sufficient to ensure satisfactory transient and steady-state performance. On the contrary, its inclusion often causes an undesired initial jump in the control signal due to the proportional action on the derivative of the error.
However, for systems with significant time delay or higher-order dynamics, the derivative component may become beneficial. In such cases, it can partially compensate the phase lag introduced by the delay and help reduce overshoot or oscillatory behavior in the closed-loop response. Nevertheless, it is not strictly necessary for stability if appropriate tuning and filtering are applied.

\section*{Conclusions and Discussion}

This paper presented a methodology for the analysis and control of nonlinear systems based on piecewise linear approximation. 
The approach provides a mathematically tractable framework that preserves essential dynamical properties of the original nonlinear model, such as continuity and boundedness, while allowing the use of well-established linear control techniques. 
By replacing the nonlinear function with a collection of locally linear segments defined over simplices, the system becomes analytically more accessible without sacrificing fidelity in regions of practical interest.

The piecewise linearization ensures Lipschitz continuity and enables the use of transfer-function-based methods for controller synthesis. 
As the partition resolution \(h\) increases, the approximation converges uniformly to the original nonlinear dynamics, leading to improved correspondence between the linearized and true system responses. 
This property allows for effective PID controller design using optimization-based methods such as Particle Swarm Optimization (PSO), while maintaining the nonlinear character of the closed-loop behavior. 
The numerical results confirm that the PID parameters optimized for the piecewise linear model yield satisfactory performance when applied to the original nonlinear system, achieving small steady-state error and limited overshoot.

An additional advantage of the piecewise linearization framework lies in its analytical transparency. 
For the linearized subsystems, quantities such as the \emph{initial jump} and \emph{initial slope} of the system response at \(t = 0^+\) can be computed explicitly, following the procedure derived in Remark~\ref{initialJumpandSlope} and Appendix~\ref{app:IJandIS}. 
These parameters are of practical importance for control implementation, as they provide valuable information about transient behavior and actuator demands immediately after control activation. 
Such explicit characterization would be analytically intractable for a fully nonlinear model, underscoring the practical utility of the proposed approximation.

A practical trade-off remains between model accuracy and computational complexity, since increasing \(h\) enhances precision at the cost of higher computational demands. Nevertheless, the proposed framework offers an analytically grounded and computationally efficient alternative to direct nonlinear optimization, providing both theoretical insight and practical applicability.

Beyond the PSO-based optimization adopted in this work, other metaheuristic and nature-inspired algorithms---such as hybrid PSO-GWO (Grey Wolf Optimizer) schemes or evolutionary-swarm com\-bi\-na\-tions---have also been reported in recent literature \cite{Siddiqi2025}. These approaches incorporate adaptive exploration-exploitation strategies to enhance convergence in high-dimensional or strongly nonconvex problems. However, for the considered piecewise linear control formulation, the optimization landscape is relatively smooth and low-dimensional, making PSO sufficiently robust and computationally efficient. 

Future research may extend this methodology to more complex nonlinear systems or control architectures, where richer optimization dynamics and parameter coupling could further benefit from hybrid or adaptive strategies.

\appendix 
\section{Proof of Lemma~\ref{Lipschitz1}}\label{app:Lipschitz}    
To prove that \( \tilde{f} \) is globally Lipschitz, we show that  
\[
|\tilde{f}(x) - \tilde{f}(y)| \leq \tilde{L} |x - y| \quad \text{for all } x, y \in \mathbb{R}^n,
\]  
where \( |\cdot| \) denotes the Euclidean norm.  

Let \( x = [x_1, x_2, \dots, x_n]^\top \) and \( y = [y_1, y_2, \dots, y_n]^\top \). Then,  
\[
\tilde{f}(x) - \tilde{f}(y) = \begin{bmatrix} x_2 - y_2 \\ x_3 - y_3 \\ \vdots \\ x_n - y_n \\ f(y) - f(x) \end{bmatrix}.
\]  
The Euclidean norm of this difference is  
\[
|\tilde{f}(x) - \tilde{f}(y)| = \sqrt{\sum_{i=2}^n (x_i - y_i)^2 + (f(y) - f(x))^2}.
\]
The Euclidean norm of the input difference is  
\[
|x - y| = \sqrt{\sum_{i=1}^n (x_i - y_i)^2}.
\]  
Using the definition of Lipschitz continuity for \( f \), we have  
\(
|f(x) - f(y)| \leq L |x - y|. 
\)
Substituting the Lipschitz bound for \( f \),  
\[
|\tilde{f}(x) - \tilde{f}(y)| \leq \sqrt{\sum_{i=2}^n (x_i - y_i)^2 + L^2 |x - y|^2}.
\]
Separating \( |x - y|^2 \) into its components 
\[
|x - y|^2 = (x_1 - y_1)^2 + \sum_{i=2}^n (x_i - y_i)^2.
\]  
Therefore,  
\[
|\tilde{f}(x) - \tilde{f}(y)| \leq \sqrt{|x - y|^2 - (x_1 - y_1)^2 + L^2 |x - y|^2}.
\]
Factoring \( |x - y|^2 \) yields
\[
|\tilde{f}(x) - \tilde{f}(y)| \leq \sqrt{|x - y|^2 \big(1 - \frac{(x_1 - y_1)^2}{|x - y|^2} + L^2\big)}\leq \left(\sqrt{1 - \frac{(x_1 - y_1)^2}{|x - y|^2} + L^2}\,\right)|x - y|.
\]
Maximizing the coefficient inside the square root gives the Lipschitz constant for \( \tilde{f} \),  
\[
\tilde{L} = \sqrt{(n-1) + L^2}.
\]  

\section{Derivation of Initial Jump and Initial Slope for the General PID-Driven First-Order System}\label{app:IJandIS}
The following closed-loop equation is obtained from (\ref{eq:first_order}) by applying it to the first-order system (\ref{linear:control}) with parameters 
$a_h=a$ and $b_h=0$:
\begin{equation*}
	\begin{aligned}
		(1 + b K_d) y''(t) + (a + b K_p) y'(t) + b K_i y(t) 
		&= b \left(K_p \delta(t) +  K_i \eta(t) +  K_d \delta'(t)\right),
	\end{aligned}
\end{equation*}
where $\eta(t)$ denotes the Heaviside unit step function and $\delta(t)$ is the Dirac delta. 

To determine the initial jump $y(0^+)$ and the initial slope $y'(0^+)$, we represent the solution in the form
\begin{equation*}
	y(t) = \tilde y(t) \, \eta(t),
\end{equation*}
where $\tilde y(t)$ is a smooth function. Using distributional derivatives,
\begin{align*}
	y'(t) &= \tilde y'(t) \eta(t) + \tilde y(0)\, \delta(t), \\
	y''(t) &= \tilde y''(t) \eta(t) + \tilde y'(0)\, \delta(t) + \tilde y(0)\, \delta'(t).
\end{align*}

\emph{Initial Jump.} Collecting coefficients at $\delta'(t)$ on both sides,
\[
(1 + b K_d) \tilde y(0) \, \delta'(t) = b K_d \, \delta'(t) \quad \Longrightarrow \quad
y(0^+) = \tilde y(0) = \frac{b K_d}{1 + b K_d}.
\]

\emph{Initial Slope.} Next, collect coefficients at $\delta(t)$,
\[
(1 + b K_d) \tilde y'(0) + (a + b K_p) \tilde y(0) = b K_p.
\]
Substituting the previously found value of $\tilde y(0)$,
\[
(1 + b K_d) \tilde y'(0) + (a + b K_p) \frac{b K_d}{1 + b K_d} = b K_p.
\]
Solving for $\tilde y'(0)$ yields
\[
y'(0^+) = \tilde y'(0) = \frac{b K_p (1 + b K_d) - b K_d (a + b K_p)}{(1 + b K_d)^2}.
\]
\section{Proof of Theorem~\ref{Gronwall}}\label{app:Gronwall}
Introduce the integrating factor
\[
\mu(t) = \exp\!\Big(-\int_a^t f(s)\,ds\Big),
\]
which is positive and continuously differentiable on \([a,b]\). Multiplying the differential inequality by \(\mu(t)\) yields
\[
\mu(t)\,w'(t) \le \mu(t) f(t) w(t) + \mu(t) g(t).
\]
Since \(\mu'(t) = -f(t)\mu(t)\), the left-hand side can be rewritten using the product rule as
\[
\frac{d}{dt}\big(\mu(t) w(t)\big)= \mu(t) w'(t) + \mu'(t) w(t)= \mu(t) w'(t) - f(t)\mu(t) w(t).
\]
Combining this with the multiplied inequality gives
\[
\frac{d}{dt}\big(\mu(t) w(t)\big) \le \mu(t) g(t).
\]
Integrating both sides from \(a\) to \(t\) yields
\[
\mu(t) w(t) - \mu(a) w(a) \le \int_a^t \mu(s) g(s)\,ds.
\]
Since \(\mu(a)=\exp(0)=1\), we obtain
\[
\mu(t) w(t) \le w(a) + \int_a^t \mu(s) g(s)\,ds.
\]
Multiplying by \(\mu(t)^{-1}=\exp\!\big(\int_a^t f(s)\,ds\big)\) recovers \(w(t)\):
\[
w(t) \le w(a)\exp\!\Big(\int_a^t f(s)\,ds\Big)
+ \int_a^t g(s)\exp\!\Big(\int_s^t f(\tau)\,d\tau\Big)\,ds.
\]
This establishes the claimed bound.

\end{document}